\documentclass[12pt]{article}
\usepackage{latexsym, amsfonts, amsmath}

\usepackage[active]{srcltx}
\usepackage{epsfig}
\usepackage{graphicx}

\newtheorem{thm}{Theorem}[section]

\newtheorem{lem}[thm]{Lemma}

\newtheorem{rem}[thm]{Remark}

\def\convd{\stackrel{\cal D}{\rightarrow}}
\def\convp{\stackrel{P}{\rightarrow}}

\def\ex{{\rm E\,}}
\def\I{{\rm I\,}}

\def\var{\mathop{\rm Var}\nolimits}

\def\cov{\mathop{\rm Cov}\nolimits}

\def\bias{\mathop{\rm bias}\nolimits}

\begin{document}

\bibliographystyle{plain}

\title {Nonparametric Kernel Density Estimation\\
for Univariate Current Status Data}

\author {\begin{tabular}{c @{\hspace{1cm}} c }
Bert van Es & Catharina Elisabeth Graafland \\
{\normalsize Korteweg-de Vries Institute for Mathematics} & {\normalsize Meteorology Group}\\
{\normalsize University of Amsterdam} & {\normalsize Institute of Physics of Cantabria (CSIC-UC)}\\
{\normalsize Science Park 105-107,
 1098 XG Amsterdam,} & {\normalsize Avenida de los Castros s/n,}\\
{\normalsize P.O. Box  94248, 1090 GE Amsterdam} & {\normalsize 39005 Santander}\\
{\normalsize The Netherlands} &  {\normalsize Spain} \\
\end{tabular}
 }


\maketitle

\begin{abstract}
We derive estimators of the density of the event times of current status data. The estimators are derived for the situations where the distribution of the observation times is known  and where this distribution is unknown. The density estimators are constructed from kernel estimators of  the density of transformed current status data, which   have a distribution similar to uniform deconvolution data. Expansions of the expectation and variance as well as asymptotic normality are derived. A reference density based bandwidth selection method is proposed. A simulated example is presented.
\\[.5cm]

\end{abstract}

\section{Introduction}

The univariate current status (UCSD) problem, or type I interval censoring problem, can be formulated in the following way. Let $X_1,\ldots,X_n$ be unobservable variables of interest. Let $T_1,\ldots T_n$ be i.i.d. random variables with known density. We assume that $T_i$ is independent of $X_i$. For $T_1,\ldots, T_n$ it is known whether the unobservable variables $X_1,\ldots,X_n$ are smaller or larger than the corresponding $T_i$. In other words, writing $\Delta_i=\I_{[X_i\leq T_i]}$, the current status data $(T_1,\Delta_1),\ldots,(T_n,\Delta_n) $ are observed. Our aim is to estimate the probability density $f$ of the unobserved $X_i$.\\
An interpretation of the random variables $T_i$ as random time instants, the observation times, and the $X_i$ the event times,  explains   the term current status. The examples below show that the variables $T_1,\ldots T_n$ are not restricted to be time instants.\\

To illustrate their widespread applicability we give three examples of UCSD problems.

Given a population of children one can be interested in the age of weaning. See for instance Diamond and McDonald (1991) and  Grummer-Strawn (1993). Weaning is the process of gradually introducing an infant to its adult diet and withdrawing the supply of its mothers milk. In this problem $X_i$ is the age of an individual child when it starts the process of weaning. $T_i$ is the age of an individual child in this population. $\Delta_i$ is defined as $\Delta_i=\I_{[X_i\leq T_i]}$. Hence $\Delta_i= 1$ if the child is already weaning at observation and $\Delta_i= 0$ if it takes nothing but its mothers milk.
This problem is worth to be analysed by current status data methods, not only when the data of $X_i$ are missing, but also when the data are available, because of the inaccurate character of the initial data.

A second problem using USCD, formulated by Milton Friedman in Stat. Research Group (1947), concerns proximity fuzes. A proximity fuze is a device that detonates a munitions explosive material automatically when the distance to a target becomes smaller than a predetermined value. One can be interested in the maximum distance over which a proximity fuze operates.
In this problem $X_i$ is the maximum distance over which an individual proximity fuze operates and $T_i$ is the nearest distance the proximity fuze reaches with respect to its target. $\Delta_i$ is defined as $\Delta_i=\I_{[X_i\leq T_i]}$. Hence $\Delta_i= 1$ if the fuze operates at a distance above the nearest distance relative to its target and $\Delta_i= 0$ if the fuze operates only at a distance that is smaller than the nearest distance it has reached relative to its target.

The third and last problem is about estimation of the age of incidence of a non fatal human disease such as Hepatitis A from a cross-sectional sample. This problem was described in Keiding (1991). In this problem $X_i$ is the age of incidence of an individual by a non fatal human disease. $T_i$ is the age at which a diagnostic test is carried on the individual.

Recall $\Delta_i=\I_{[X_i\leq T_i]}$. Hence $\Delta_i= 1$ if the individual was already ill at the moment of the diagnostic test and $\Delta_i= 0$ if the individual gets ill later. Here we assume that the members of our population are all infected during their lifetime and that this infection causes illness. Data can be drawn from a population of recent and former patients. The age at which the first diagnostic test is performed can be taken as data points   $T_i$.  \\

For the problem of estimation of the density of the $X_i$, estimators have been proposed by Groeneboom, Witte and Jongbloed (2010) and Witte (2011). These estimators are based on the nonparametric maximum likelihood estimator (NPMLE) of the distribution function of the $X_i$. This estimator is subsequently smoothed by kernel based techniques to obtain  density estimators. For estimation for current status linear regression models see Groeneboom and Hendrickx (2016).

In this paper a new density estimator for the USCD problem, based on an inversion technique similar to that in Van Es (2011) for uniform deconvolution, will be constructed. This density estimator will turn out to have similar asymptotic properties as the smoothed NPMLE estimator.
An advantage of our approach is the possibility to expand the theory to bi and multivariate current status data more naturally. Extending the kernel smoothed univariate maximum likelihood estimator of the density  to the bivariate current status context is more involved. In Groeneboom (2013) and Section 12.3 of Groeneboom and Jongbloed (2014)   estimators for the bivariate distribution function are proposed. Maathuis (2005) presents an algorithm for computation of the NPMLE in the bivariate model.\\

The paper  is organized as follows. In Section \ref{const} we introduce a transformation of the UCSD date to random variables $V_i$ which have a distribution similar to uniform deconvolution data. For this type of data two inversion formulas, expressing the density  of the unobserved $X_i$ in terms of the density of the observed $V_i$, can be derived. By plugging in a kernel estimator for the density of the $V_i$ we then obtain two estimators of the density of the $X_i$. These are then combined in a convex combination   with estimated weights, minimizing the asymptotic variance, to obtain the final density estimator. Up to here we have assumed that the density of the $T_i$ is known. The next step is to estimate this density to plug in its estimator in the previous one. Thus, in Section \ref{unknown}, we get a final estimator for the case that the density of the $T_i$ is unknown. For these estimators expansions of the bias and variance, where possible, and asymptotic normality are derived. As an illustration we present a simulated example in Section \ref{simulation}. The more technically involved proofs are postponed to Section \ref{proofs}.

\section{Construction and results}\label{const}

\subsection{A transformation of the univariate current status data}\label{ParUCSD}

The basic step in our approach is a transformation of the current status data to data of which the distribution is similar to uniform deconvolution data.
Consider the following transformation of the points $T_i$,
\begin{equation}V_i=
\left\{
\begin{array}{ll}
T_i+1 &,\ \mbox{if}\ \Delta_i=0,\\
T_i &,\ \mbox{if}\ \Delta_i=1.
\end{array}
\right.
\end{equation}
The next lemma derives the distribution of the transformed data.
\begin{lem}\label{lemtransformation}
Assume that the distribution of the variables $X_i$, with distribution function $F$, is concentrated on $[0,1]$.
Assume that the $T_i$ have a density $q$.
If $T_i$ is supported on $[0,1]$, then the $V_i$ have a density $g$, given by
\begin{equation}\label{gdens1}
g(v)={\bar q}(v)\Big(F(v)-F(v-1)\Big), \quad v \in [0,2],
\end{equation}
with the function $\bar q$ defined by
\begin{equation}
{\bar q}(v )=q(v)+q(v-1).
\end{equation}
\end{lem}

\noindent{\bf Proof}\\
Note that the probabilities for $\Delta_i$, given the value of $T_i$, are given by
\begin{align}
&P(\Delta_i=0|T_i= t_i)=1-F(t_i),\nonumber\\
&P(\Delta_i=1|T_i= t_i)=F (t_i).\label{deltas}
\end{align}
We omit the subscript $i$. For $v\in [0,1] $ we have
\begin{align*}
P(V \leq v ) &=\sum_{k=0}^1 P(V \leq v , \Delta=k)
 =P(V\leq v , \Delta=1)\\
 &=\int_0^1 P(V \leq v , \Delta=1| T =t)q(t )dt
 =\int_0^1 P(T\leq v , \Delta=1| T =t)q(t )dt\\
 &=\int_0^{v} P(\Delta=1|T=t)q(t )dt
  =\int_0^{v } F(t )q(t )dt .
\end{align*}
By the support restriction on
the distribution induced by $F$ and $q$ this confirms (\ref{gdens1}) for the given values of $v$.

Let us also check the claim on the interval $[1,2]$.
For $u\in [0,1] $ we have
\begin{align*}
P(1\leq V \leq& 1+u )  = P(1\leq V \leq 1+u, \Delta=0)+P(1\leq V \leq 1+u , \Delta=1)\\
&=P(1\leq V \leq 1+u , \Delta=0)
 =\int_0^1 P(1\leq V \leq 1+u , \Delta=0| T =t)q(t )dt \\
&=\int_0^1 P(T\leq u, \Delta=0| T =t)q(t)dt
 =\int_0^{u} P(\Delta=0|T=t)q(t )dt\\
&=\int_0^{u } (1-F(t ))q(t )dt
 =\int_1^{1+u } (1-F(t-1 ))q(t-1 )dt.
\end{align*}
Again by the support restriction on
the distribution induced by $F$ and $q$ this confirms (\ref{gdens1}) for the values of $v$ in $[1,2]$.$\hfill\Box$\\

This lemma reveals a connection between the UCSD problem and the uniform deconvolution problem. In the uniform deconvolution problem we have observations with density
\begin{equation}\label{uddensity}
g(v)= F(v)-F(v-1).
\end{equation}
The USCD problem is equivalent to the uniform deconvolution problem if the variables $T_i$ have a uniform distribution on $[0,1]$, since then $q \equiv 1$ on $[0,1]$. In this case $V_i$ is in distribution equal to $X_i + Z_i$, with $Z_i \sim Un[0,1)$.

\subsection{Inversion formulas}\label{ParInvForm}
We deduce the following inversion formulas from Lemma \ref{lemtransformation}. These formulas express the density $f$ and distribution function $F$ in terms of the density $g$ of the transformed  current status data.

\begin{lem}\label{inversionlemma}
If $g$ is of the form (\ref{gdens1}) and $q$ is strictly positive on $[0,1]$, then we have for $x \in [0,1]$

\begin{equation}\label{inversion1}
F(x) = \frac{g(x)}{q(x)},
\end{equation}
\begin{equation}\label{inversion2}
F(x) = 1 - \frac{g(x+1)}{q(x)}.
\end{equation}
Furthermore, if we assume that $F$ and $q$ are differentiable, we have for $0 < x < 1$
\begin{equation}\label{inversion3}
\quad f(x) = \frac{1}{q(x)}g'(x) - \frac{q'(x)}{q^{2}(x)}g(x),
\end{equation}
\begin{equation}\label{inversion4}
\quad f(x) =  - \left( \frac{1}{q(x)}g'(x+1) - \frac{q'(x)}{q^{2}(x)}g(x+1) \right).
\end{equation}
\end{lem}

\noindent{\bf Proof}\\
Note that $v\in [0,1] $ implies $F(v-1) = 0$ and $q(v-1) = 0$. Hence for $v\in [0,1] $ equation (\ref{gdens1}) turns into $g(v) = F(v)q(v)$. Rewriting this gives the first inversion formula (\ref{inversion1}).\\
Also $v\in [1,2] $ implies $F(v) = 1$ and $q(v) = 0$. Hence for $v\in [1,2] $ equation (\ref{gdens1}) turns into $g(v) = (1-F(v-1))q(v-1)$. Rewriting this gives $F(v-1) = 1 - \frac{g(v-1)}{q(v-1)}$. Substituting $x$ for $v-1$ gives the second inversion formula (\ref{inversion2}) for $x \in [0,1]$.\\
Differentiating these formulas yields the two formulas (\ref{inversion3}) and (\ref{inversion4}).\hfill$\Box$\\

From now on we assume the density $q$ to be differentiable and strictly positive on $[0,1]$. We also assume that the distribution function $F$ has a density $f$.

\subsection{Two estimators of the density function}\label{ParConstrEstDen}
We start with the construction of two different estimators using the inversion formulas in Lemma \ref{inversionlemma}. Our aim is to combine these estimators by a convex combination to get an `optimal' estimator for the density function $f$ of the unobservable variables of interest $X_1,\ldots,X_n$.

Note that the inversion formulas in Lemma \ref{inversionlemma} yield equal $F(x)$ and $f(x)$ if $g$ is exactly of the form (\ref{gdens1}). However for arbitrary $g$, for example an estimator of $g$, which is not of the form (\ref{gdens1}), the inversions will in general not coincide. Also they may not yield distribution functions or densities.\\

We now get two different estimators of $f(x)$ from (\ref{inversion3}) and
(\ref{inversion4}) in the following way. For $g(x)$ we substitute the kernel density
estimator, given by
\begin{equation}\label{kerneldensityestimator}
g_{nh}(x)=\frac{1}{n}\,\sum_{i=1}^n \frac{1}{h}\, w\Big(\frac{x-V_i}{h}\Big).
\end{equation}
Here $h>0$ is called the bandwidth which controls the roughness of the estimate and the function $w$ is called the kernel function. We impose the following condition on the kernel function $w$.

\medskip

\noindent{\bf Condition W}

The function $w$ is a continuously differentiable symmetric
probability density function with support [-1,1].\\

\noindent General books on kernel estimation are for instance Prakasa Rao (1983), Silverman (1986) and Wand and Jones (1995).

For $g'(x)$ we substitute the derivative of this kernel estimator, $g_{nh}'(x)$,
\begin{equation}\label{kerneldensityderivativeestimator}
g_{nh}'(x)=\frac{1}{n}\,\sum_{i=1}^n \frac{1}{h^2}\, w'\Big(\frac{x-V_i}{h}\Big).
\end{equation}

In this way we derive two different estimators of $f(x)$ from (\ref{inversion3}) and
(\ref{inversion4}).
We call the estimators respectively the left and right estimator, written as $f_{nh}^{-}$ and $f_{nh}^{+}$. We have
\begin{equation}\label{estimator1}
\quad f^{-}_{nh}(x) = \frac{1}{q(x)} g_{nh}'(x) - \frac{q'(x)}{q^{2}(x)}g_{nh}(x)
\end{equation}
and
\begin{equation}\label{estimator2}
\quad f^{+}_{nh}(x) = - \left( \frac{1}{q(x)}g_{nh}'(x+1) - \frac{q'(x)}{q^{2}(x)}g_{nh}(x+1) \right).
\end{equation}

Note that these estimators coincide with the estimators in the uniform deconvolution model as obtained in Van Es (2011) when the $T_i$ are taken to be uniform, i.e. $q(x) \equiv 1$ on $[0,1]$.

In the next section we derive  expansions of the expectation and variance of $f^{-}_{nh}(x)$ and $f^{+}_{nh}(x)$, showing that they are  both consistent estimators of $f(x)$.

\subsubsection{Expectation and variance}\label{ParExpVarleftright}

The next two theorems give the expansions for the left and right estimator. The proofs can be found  in Section \ref{proofs}.

\begin{thm}\label{exvarfnhmin}
Assume that the kernel function $w$ satisfies Condition W, that the density $f$ is twice continuously differentiable, and that the density $q$ is three times continuously differentiable, then we have for $0 < x < 1$
\begin{equation}\label{exfnhmin}
\ex f^{-}_{nh}(x) = f(x) + \frac{1}{2}h^2\int v^2w(v)dv \: b^-(x) + o(h^2),
\end{equation}
and
\begin{equation}\label{varfnhmin}
\var f^{-}_{nh}(x) = \frac{1}{q(x)}\frac{1}{nh^3}F(x)\int w'(u)^2du + o\left(\frac{1}{nh^3}\right),
\end{equation}
where
\begin{equation}\label{bmin}
b^-(x) = \frac{1}{q(x)}g'''(x) - \frac{q'(x)}{q^{2}(x)}g''(x).
\end{equation}
\end{thm}

\begin{thm}\label{exvarfnhplus}
Assume that the kernel function $w$ satisfies Condition W, that the density $f$ is twice continuously differentiable, and that the density $q$ is three times continuously differentiable, then we have for $0 < x < 1$
\begin{equation}\label{exfnhplus}
\ex f^{+}_{nh}(x) = f(x) + \frac{1}{2}h^2\int v^2w(v)dv \: b^+(x) + o(h^2),
\end{equation}
and
\begin{equation}\label{varfnhplus}
\var f^{+}_{nh}(x) = \frac{1}{q(x)}\frac{1}{nh^3}(1-F(x))\int w'(u)^2du + o\left(\frac{1}{nh^3}\right),
\end{equation}
where
\begin{equation}\label{bplus}
b^{+}(x) = -\frac{1}{q(x)}g'''(x+1) + \frac{q'(x)}{q^{2}(x)}g''(x+1).
\end{equation}
\end{thm}

These theorems show that both estimators have a bias of order $h^2$ and a variance of order $1/(nh^3)$. If we minimize the mean squared error with respect to $h$ we get an optimal rate $n^{-1/7}$ for the bandwidth and an optimal mean squared error of order $n^{-4/7}$. These rates will also appear in the asymptotics of our later estimators which are essentially convex combinations of the left and right estimator. The results also show that if we choose the bandwidth suboptimal, i.e. $h \ll n^{-1/7}$ then the bias is negligible compared to the standard deviation.

\subsection{Convex combination of the left and right estimator}\label{ParCombleftright}
The results in the previous section show that the left estimator has a smaller variance than the right estimator for small values of $x$, because of the factors $F(x)$ and $1-F(x)$. For large values of $x$ the converse is true. Hence it makes sense to construct a convex linear combination of the two estimators. We define the combined estimator by
$$
f^{t}_{nh}(x) = t f^{-}_{nh}(x) + (1-t)f^{+}_{nh}(x),
$$
for some fixed $0\leq t\leq 1$. This factor $t$ is later chosen dependent on $x$ to dampen the effect of larger variances at the endpoints.\\

One could wonder if the bias of this combined estimator will not be bigger than the bias of $f^{-}_{nh}(x)$ and/or $f^{+}_{nh}(x)$. By linearity of the bias, the bias of the new estimator for $t \in (0,1)$ will lie somewhere between the bias of $f^{-}_{nh}(x)$ and $f^{+}_{nh}(x)$. The exact value depends on $t$. In Section \ref{ConsEstF} we will find an optimal value for $t$, which we call $t^*$. The value $t^*$ is optimal in the sense that it minimizes the mean squared error of $f^{t}_{nh}(x)$ with respect to $t$. It is possible that $f^{t^*}_{nh}(x)$ will have a larger bias than $f^{-}_{nh}(x)$ or $f^{+}_{nh}(x)$, but it will have lower or equal mean squared error for certain.

\subsubsection{Expectation and  variance}\label{ParExpVarComb}
We will derive expansions of the expectation and the variance of the combined estimator as well as its asymptotic normality. The next theorem states the expansions of the expectation and variance of $f^{t}_{nh}(x)$.
For the proof we refer to Section 6.
\begin{thm}\label{ExVarfnht}
Assume that the kernel function $w$ satisfies Condition W, that the density $f$ is twice continuously differentiable, and that the density $q$ is three times continuously differentiable, then we have for $0 < x < 1$
\begin{equation}\label{exfnht}
\ex f^{t}_{nh}(x) = f(x) + \frac{1}{2}h^2\int v^2w(v)dv \left(t b^-(x) + (1-t)b^+(x)\right) +o(h^2),
\end{equation}
and
\begin{equation}\label{varfnht}
\var f^{t}_{nh}(x) = \frac{1}{q(x)}\frac{1}{nh^3}\left(t^2F(x) + (1-t)^2(1-F(x)\right)\int w'(u)^2du + o\left(\frac{1}{nh^3}\right),
\end{equation}
where the functions $b^-$ and $b^+$ are defined by (\ref{bmin}) and (\ref{bplus}).
\end{thm}

\begin{rem}{\rm
Note that for uniformly distributed variables $T_i$ on $[0,1]$, i.e. $q \equiv 1$ on $[0,1]$, we have $b^-(x) = b^+(x) = f''(x)$. In this case the expectation of $f^{t}_{nh}$ does not depend on $t$.}
\end{rem}

\subsubsection{Asymptotic normality}\label{asknownq}
We now derive   asymptotic normality of the combined estimator. Recall that, for some fixed $0\leq t\leq 1$, we have
$$
f^{t}_{nh}(x) = t f^{-}_{nh}(x) + (1-t)f^{+}_{nh}(x)
$$
and note that asymptotic normality for $f^{t}_{nh}(x)$ implies  asymptotic normality of our estimators $f^{-}_{nh}(x)$ and $f^{+}_{nh}(x)$ by taking $t$ equal to zero and one.

The proof is given in Section \ref{proofs}.

\begin{thm}\label{AsNormfnht}
Assume that Condition $W$ is satisfied and that $f$ is bounded on a
neighbourhood of $x$. Then, as $n\to\infty$, $h\to 0$ and $nh^3 \to \infty$, we have for $0 < x < 1$
$$
\sqrt{nh^3} (f^{t}_{nh}(x)-\ex f^{t}_{nh}(x))\convd
N(0,\sigma_t^2),
$$
with
\begin{equation}\label{asvar}
\sigma_t^2 = \frac{1}{q(x)}\left(t^2F(x) + (1-t)^2(1-F(x)\right)\int w'(u)^2du.
\end{equation}
\end{thm}

\bigskip
The asymptotic bias and variance of $f^{t}_{nh}(x)$ depend on the factor $t$. Below we will derive an optimal choice for $t$.

\subsection{The final estimator with $q$ known}\label{ConsEstF}

The results for the combined estimator $f^{t}_{nh}(x)$ show that the bias of the estimator is quite complicated in the sense that it depends on the unknown second and third derivative of $g$. In the case that $q$ is unknown the bias additionally depends on $q$ and its derivative. On the other hand the asymptotic   variance (\ref{asvar})
is relatively simple in its dependence on $F$ and $q$.

Minimizing the asymptotic variance with respect to $t$ we get the optimal value
$ t^* = 1-F(x)$.
For this choice of $t$ the asymptotic mean squared error has the optimal rate $n^{-4/7}$ if $h$ has the optimal order $n^{-1/7}$. It does not minimize the asymptotic mean squared error. Minimizing the mean squared error would yield an unpractical dependence on $g$ and $q$. For sub optimal $h\ll n^{-1/7} $ however we have an optimal mean squared error, which in this case equals the asymptotic variance.

 Since the optimal $t^*$ still depends on unknown quantities we have to plug in  an estimator. Hence our final estimator will be equal to  an $f^t_{nh}(x)$ with $t = 1-\hat{F}_{n}(x)$, where $\hat{F}_{n}(x)$ is an estimator for $F(x)$ that is consistent in mean squared error. Note that $1-\hat{F}_n(x)$ is not generally in $[0,1]$ since this is not required of $\hat{F}_n(x)$. Hence we now have constructed our final estimator. It is given by

\begin{equation}\label{finalest}
f_{nh}(x) = (1-\hat{F}_{n}(x)) f^{-}_{nh}(x) + \hat{F}_{n}(x)f^{+}_{nh}(x).
\end{equation}
Note that we use the fact that the density $q$ is known since it appears in the construction of the left and right estimator. Later on we will drop this assumption and  present an estimator for the case where this density is not known.

\subsubsection{Asymptotic normality}

Our main theorem for the situation where $q$ is known gives an expansion for the bias of $f_{nh}(x)$ and it establishes asymptotic normality. Its proof is postponed to Section \ref{proofs}.
\begin{thm}\label{theorem2}
Assume that Condition $W$ is satisfied, that $f$ is twice continuously differentiable, $q$ is three times continuously differentiable on a neighbourhood
of $x$, and that $\hat F_n(x)$ is an estimator of
$F(x)$ with
\begin{equation}\label{fhatcond}
\ex (\hat F_n(x)-F(x))^2\to 0.
\end{equation}
Then, as $n\to\infty$ and $h = O(n^{-1/7})$, we have $nh^3 \to\infty$, $h\to 0$ and we have for $0<x<1$
\begin{equation}
\sqrt{nh^3}(f_{nh}(x)-\ex f_{nh}(x))\convd N(0,\sigma^2),
\end{equation}
with
\begin{equation}
\sigma^2=F(x)(1-F(x))\frac{1}{q(x)}\int w'(u)^2du.
\end{equation}
Furthermore, if
\begin{equation}\label{fhatcond2}
\ex (\hat F_n(x)-F(x))^2=o(nh^7)
\end{equation}
then
\begin{equation}\label{biasfinal}
\ex  f_{nh}(x) = f(x) + \frac{1}{2}h^2\int v^2w(v)dv\left((1-F(x))b^-(x) + F(x)b^+(x)\right)+ o(h^2),
\end{equation}
where the functions $b^-$ and $b^+$ are defined by (\ref{bmin}) and (\ref{bplus}).

\noindent  If $h \ll n^{-1/7}$ then $f_{nh}(x)$ is t-optimal and we have
\begin{equation}
\sqrt{nh^3}(f_{nh}(x)- f(x))\convd N(0,\sigma^2),
\end{equation}
as $n\to\infty$.
\end{thm}

\begin{rem}{\rm
The bias term in (\ref{biasfinal}) depends indirectly on $q$ and $f$. After some computation we get
\begin{align*}
b^-(x)&=\frac{1}{q(x)^2}\Big(q(x)q'''(x)F(x)+3q(x)q''(x)f(x)+2q(x)q'(x)f'(x)+q^2(x)f''(x)\\
&\quad\quad\quad-q'(x)q''(x)F(x)
-2q'(x)^2f(x)\Big),\\
b^+(x)&=\frac{1}{q(x)^2}\Big(-q(x)q'''(x)(1-F(x))+3q(x)q''(x)f(x)+2q(x)q'(x)f'(x)+q^2(x)f''(x)\\
&\quad\quad\quad+q'(x)q''(x)(1-F(x))-2q'(x)^2f(x)\Big).
\end{align*}
Hence
\begin{align}
(1-F(x))&b^-(x) + F(x)b^+(x)\nonumber\\
&=\frac{1}{q(x)^2}\Big(3q(x)q''(x)f(x)+2q(x)q'(x)f'(x)+q^2(x)f''(x) -2q'(x)^2f(x)  \Big).\label{biasreduced}
\end{align}
For the relatively simple case where $q$ is the uniform density on [0,1] this expression equals $f''(x)$.}
\end{rem}

\subsubsection{Consistent estimator of the distribution function}\label{chaconsistentF}
Theorem \ref{theorem2} is based on the assumption that an estimator that is consistent in mean squared error exists for the distribution function $F$. In this section we show that such an estimator can be found easily. In fact there are many suitable estimators of $F$.\\

The construction of the estimator starts with the inversion formulas (\ref{inversion1}) and (\ref{inversion2}) obtained in Section \ref{ParInvForm}. For $g$ we substitute again the kernel estimator $g_{nh}$, defined in (\ref{kerneldensityestimator}), as we did in the construction of $f^-_{nh}$ and $f^+_{nh}$. This method leads to two different estimators for $F$ for $x \in [0,1]$,

\begin{equation}\label{Festimators}
F_{nh}^{-}(x) = \frac{g_{nh}(x)}{q(x)} \quad
\mbox{and}\quad F_{nh}^{+}(x) = 1 - \frac{g_{nh}(x+1)}{q(x)}.
\end{equation}
We define the estimator $F^t_{nh}(x)$ for $x \in [0,1]$ as
\begin{equation}
F^t_{nh}(x) = tF^-_{nh}(x) + (1-t)F^+_{nh}(x)
\end{equation}
for some fixed $0\leq t\leq 1$. In the sequel we show that the estimator $F^t_{nh}$ is consistent in mean squared error. This is proven with the help of expansions of the expectation and the variance of $F^t_{nh}$. These expansions are stated in the next theorem of which the proof can be found in Section \ref{proofs}.\\

\begin{thm}\label{Theorem3} Under the assumptions of Theorem \ref{theorem2}, that is, if Condition $W$ is satisfied, $f$ is twice continuously differentiable and $q$ is three times continuously differentiable on a neighbourhood
of $x$, we have
\begin{equation}\label{a8}
\ex F_{nh}^{t}(x) = F(x) + \frac{1}{q(x)}\frac{1}{2}h^2\left(tg''(x)-(1-t)g''(x+1)\right)\int v^2w(v)dv + o(h^2)
\end{equation}
and
\begin{equation}\label{a9}
\var F^{t}_{nh}(x) = \frac{1}{q(x)}\frac{1}{nh}(t^2 F(x) +(1-t)^2(1-F(x))\int
w(v)^2dv+ o\Big( \frac{1}{nh}\Big).
\end{equation}
\end{thm}
\bigskip
With the help of Theorem \ref{Theorem3} we are able to check that the estimator $F^{t}_{nh}$ is consistent in MSE. We have
\begin{equation}\label{OrderFnht}
\begin{split}
\ex \left(F_{nh}^{t}(x) - F(x)\right)^2 &= \bias^2 F_{nh}^{t}(x) + \var F^{t}_{nh}(x)\\
&= \left(\frac{1}{q(x)}\frac{1}{2}h^2\left(tg''(x)-(1-t)g''(x)\right)\int v^2w(v)dv + o(h^2)\right)^2\\
&\quad + \frac{1}{q(x)}\frac{1}{nh}(t^2 F(x) +(1-t)^2(1-F(x))\int
w(v)^2dv+ o\Big( \frac{1}{nh}\Big)\\
&= O\big(h^4\big) + O\big(\frac{1}{nh}\big)  \to 0,
\end{split}
\end{equation}
as $h \to 0$ and $nh \to \infty$.\\

The second statement in Theorem \ref{theorem2} is based on the assumption that the estimator $F_{n}(x)$ satisfies
$$\ex \left(\hat F_{n}(x) - F(x)\right)^2 = o(nh^7).$$ This assumption is stated in equation (\ref{fhatcond2}).
Note that we can choose different bandwidths $h$ in $\hat F_{n}(x)$ and $f_{nh}(x)$ as long as we meet the requirements in Theorem \ref{theorem2}. In the sequel we will write $h_{1}$ for the bandwidth in $f_{nh}(x)$ and $h_{2}$ for the bandwidth in $\hat F_{n}(x)$.  We show below in which way the estimator $F_{nh}^{t}(x)$ is able to satisfy equation (\ref{fhatcond2}) for respectively an optimal and a sub optimal choice of $h_{1}$. \\

An optimal choice of $h_{1}$, i.e. $h_{1} = cn^{-\frac{1}{7}}$ with $c>0$, reduces the assumption for $F_{nh}^{t}(x)$ to
$$\ex \left(F_{nh}^{t}(x) - F(x)\right)^2 = o(1).$$
By equation (\ref{OrderFnht}) this requirement is met by the estimator $F_{nh}^{t}(x)$ for all $h_{2}$ such that $h_{2} \to 0$ and $nh_{2} \to \infty$ as $n \to \infty$.
In particular the requirement is met for the optimal choice of $h_{2}$ in $F_{nh}^{t}(x)$, given by $h^*_{2} = cn^{-1/5}$ with $c>0$. The bandwidth $h^*_{2}$ is optimal in the sense that it minimizes the MSE (equation (\ref{fhatcond2})) with respect to $h_{2}$.\\

A sub optimal choice of $h_{1}$ requires the MSE of $F_{nh}^{t}(x)$ to be of smaller order.
It is possible to keep the optimal choice $h^*_{2} \sim n^{-1/5}$ when $h_{1}$ is chosen sub optimal. However in order to satisfy equation (\ref{fhatcond2}) one should not choose $h_{1}$ too small. Note that by equation (\ref{OrderFnht})
the MSE of $F_{nh}^{t}(x)$ with $ h = h^*_{2}$ is of order $O(n^{-4/5})$. Hence to satisfy equation (\ref{fhatcond2}) one should always choose $h_{1} \gg n^{-9/35}$.
With this restriction $F_{nh}^{t}(x)$ satisfies equation (\ref{fhatcond2}) for sub optimal $h_{1}$ and optimal $h_{2}$.

\section{The final estimator with $q$ unknown}\label{unknown}

Consider the situation where the density $q$ of the observation points $  T_i$, on
the unit interval, is not known.  We can then estimate $q$ from the data $ T_1 ,\ldots,  T_n$ by a kernel estimator $q_{n{\tilde h}}$ with a bandwidth $\tilde h$,
\begin{equation}\label{kernelq0}
	q_{n\tilde h}(x)=\frac{1}{n\tilde h}\sum_{k=1}^nw\Big (\frac{x -T_k}{\tilde h} \Big).
\end{equation}
The following theorem establishes asymptotic normality. Its proof is postponed to Section \ref{proofs}.

\begin{thm}\label{ThmUnkownq}
Assume that Condition $W$ is satisfied, that $f$ is twice continuously differentiable, $q$ is three times continuously differentiable on a neighbourhood
of $x$, and that $\hat F_n(x)$ is an estimator of
$F(x)$ with
\begin{equation}\label{mseconsistency}
\ex (\hat F_n(x)-F(x))^2\to 0.
\end{equation}
If both $h$ and $\tilde h$ are of order $n^{-1/7}$, i.e. $h=cn^{-1/7}$ and $\tilde h=\tilde cn^{-1/7}$ for  $c>0$ and $\tilde c>0$, then we have as $n\to\infty$ for $0<x<1$
\begin{equation}
n^{2/7}(f_{nh\tilde h} (x ) -f(x))\stackrel{\cal D}{\to}N(\mu
,c^{-3}\sigma^2),
\end{equation}
with
\begin{equation}\label{biasunknownq}
\mu=\frac{1}{2}\,\int_{-\infty}^\infty v^2 w(v)dv\Big(c^2((1-F(x)b^-(x)+F(x)b^+(x))
 -{\tilde c}^2\,\frac{f(x)q''(x)}{q(x)}\Big),
\end{equation}
where the functions $b^-$ and $b^+$ are defined by (\ref{bmin}) and (\ref{bplus}), and
$$
\sigma^2=F(x)(1-F(x))\frac{1}{q(x)}\int w'(u)^2du.
$$
If we choose $h$ suboptimal, i.e. $n\to \infty$, $nh^3\to \infty$, $h \ll  n^{-1/7}$ and $\tilde h=\tilde cn^{-1/7}$, then we have for $0<x<1$
\begin{equation}
\sqrt{nh^3}(f_{nh\tilde h} (x )-f(x))\stackrel{\cal D}{\to}N(0,\sigma^2).
\end{equation}
\end{thm}

\begin{rem}{\rm
In the theorem we have two bandwidths $h$ and $\tilde h$ which are both of order $n^{-1/7}$.  If we substitute (\ref{biasreduced}) and choose the same bandwidths, say $h$, then the bias (\ref{biasunknownq}) reduces  to
\begin{align*}
\frac{1}{2}&h^2 \int_{-\infty}^\infty v^2 w(v)dv\Big( ((1-F(x)b^-(x)+F(x)b^+(x))
 - \frac{f(x)q''(x)}{q(x)}\Big)\\
&=\frac{1}{2}\,h^2 \int_{-\infty}^\infty v^2 w(v)dv\Big(\frac{3q(x)q''(x)f(x)+2q(x)q'(x)f'(x)+q^2(x)f''(x)-2q'(x)^2f(x)}{q(x)^2}\\
&\quad\quad\quad\quad- \frac{f(x)q''(x)}{q(x)}\Big)\\
&=\frac{1}{2}\,h^2 \int_{-\infty}^\infty v^2 w(v)dv\Big(\frac{2(x)q''(x)f(x)+2q(x)q'(x)f'(x)+q^2(x)f''(x)-2q'(x)^2f(x)}{q(x)^2}\Big).
\end{align*}
This asymptotic bias is exactly the same as the asymptotic bias of the Maximum Smoothed Likelihood density estimator of Groeneboom et al. (2010). Given that the asymptotic variance is also the same, our estimator has the same asymptotics as this estimator.  Their other estimator, the smoothed maximum likelihood estimator has the same asymptotic variance but a different bias. This bias can be smaller or larger than our bias, depending on the specific $q$ and $f$. Admittedly both their estimators yield true non negative densities while our estimator can take on negative values.}
\end{rem}

\begin{rem}{\rm
Theorem \ref{AsNormfnht}, Theorem \ref{theorem2} and Theorem \ref{ThmUnkownq} are based on the assumption that the support of the random variables $T_i$ is $[0,1]$ and the domain of the random  variables $X_i$ is $[0,1]$. These restrictions are however only given for simplicity. One can adapt the theory to the more general problem in which the variables $T_i$ have bounded support $[a,b]$ and the distribution of the variables $X_i$ is concentrated on $[a,\infty)$.
The transformation of the variables $T_i$ for $x \in [a,b]$ is as follows,
\begin{equation}\nonumber
V_i=
\left\{
\begin{array}{ll}
T_i+b-a &,\ \mbox{if}\ \Delta_i=0,\\
T_i &,\ \mbox{if}\ \Delta_i=1.
\end{array}
\right.
\end{equation} The construction of the estimators is similar as before and the theorems remain valid for $x \in [a,b]$.
}\end{rem}

\section{Bandwidth selection}
Let us first assume that the density $q$ of the observation times is known. From the properties of the estimator $f_{nh}$ stated in Theorem \ref{theorem2} we can derive the following expansion of the mean integrated squared error.
The expansion holds because the integrals are over a finite interval and since the expansions of the bias and variance still hold for converging sequences $x_n$ replacing a fixed $x$, thus rendering the expansions uniform on [0,1]. We have
\begin{align}
MISE&(h)=\ex \int_0^1(f_{nh}(x)-f(x))^2dx= \int_0^1\Big((\ex f_{nh}(x)-f(x))^2+\var f_{nh}(x)\Big) dx\nonumber\\
&=\frac{1}{4}h^4\Big(\int v^2w(v)dv\Big)^2\int_0^1\left((1-F(x))b^-(x) + F(x)b^+(x)\right)^2dx\\
&\quad+\frac{1}{nh^3}\int w'(u)^2du\int_0^1\frac{1}{q(x)}\,F(x)(1-F(x))dx+o(h^4)+o\Big(\frac{1}{nh^3}\Big).
\nonumber
\end{align}
The asymptotically optimal bandwidth, minimizing the asymptotic mean integrated squared error is given by
\begin{equation}\label{optband}
h_n^{opt}= \frac{\Big[3\Big(\int v^2w(v)dv\Big)^2\int_0^1\left((1-F(x))b^-(x) + F(x)b^+(x)\right)^2dx\Big]^{1/7}}{\Big[\int w'(u)^2du\int_0^1\frac{1}{q(x)}\,F(x)(1-F(x))dx\Big]^{1/7}}\,  n^{-1/7}.
\end{equation}
Note that this optimal bandwidth depends of the unknown density $f$ in a complicated manner.

For the relatively simple case that the observation times are uniformly distributed and $q$ is identically equal to one on [0,1], we have $g=F$ and $b^-=b^+=f^{''}$. So in this case the optimal bandwidth reduces to
\begin{equation}
h_n^{opt}= \frac{\Big[3\Big(\int v^2w(v)dv\Big)^2\int_0^1f^{''}(x)^2dx\Big]^{1/7}}{\Big[\int w'(u)^2du\int_0^1\,F(x)(1-F(x))dx\Big]^{1/7}}\,  n^{-1/7}.
\end{equation}
For general $q$ the optimal bandwidth is more involved.
The optimal bandwidth in the general case equals (\ref{optband}) with $
(1-F(x))b^-(x) + F(x)b^+(x)$ replaced by the   expression (\ref{biasreduced}) involving up to second derivatives of both $q$ and $f$.

In order to approximate the optimal bandwidth we will apply a method of reference densities, which is similar to the use of a normal reference density in direct kernel estimation as in Section 3.4.2 of Silverman (1986). Here this means that we assume a Beta$(\alpha,\beta)$ reference density for the density $f$ of $X$ and estimate its parameters. This will yield a parametric estimate of the distribution of $X$ which can be used in the optimal bandwidth (\ref{optband}).

For $\alpha>0$ and $\beta>0$ the Beta $(\alpha,\beta)$ density is given by
$$
f_{\alpha,\beta}(x)=\frac{1}{B(\alpha,\beta)}\, x^{\alpha -1}(1-x)^{\beta-1},\quad  0\leq x \leq 1,
$$
with $B(\alpha,\beta)=\int_0^1x^{\alpha -1}(1-x)^{1-\beta}$. We will write $F_{\alpha,\beta}(x)$ for the distribution function.

We use the method of moments to estimate the parameters $\alpha$ and $\beta$ from the sample of the $V_i$, using the first two moments of $X$.
Note that we have, with $V$ having the density (\ref{gdens1}),
\begin{equation}\label{mom1}
\ex \Big(\frac{V}{\bar q(V)}\Big)=\int_{-\infty}^\infty \frac{v}{\bar q(v)}\,g(v)dv=\int_{-\infty}^\infty \frac{v}{\bar q(v)}\,\bar q(v)(F(v)-F(v-1))dv
=\int_{-\infty}^\infty vg_{ud}(v)dv,
\end{equation}
where $g_{ud}$ denotes density (\ref{uddensity}) of the observations in the uniform deconvolution problem.

Similarly we get
\begin{equation}\label{mom2}
\ex \Big(\frac{V^2}{\bar q(V)}\Big)=\int_{-\infty}^\infty \frac{v^2}{\bar q(v)}\,g(v)dv=\int_{-\infty}^\infty \frac{v^2}{\bar q(v)}\,\bar q(v)(F(v)-F(v-1))dv
=\int_{-\infty}^\infty v^2g_{ud}(v)dv,
\end{equation}
From the uniform deconvolution we recall that $g_{ud}$ is the density of a random variable $X+Z$ with $X$ and $Z$ independent, $X$ having density $f$ and $Z$ equal to a Un[0,1) distributed random variable.
This gives
\begin{equation}\label{mom3}
\ex \Big(\frac{V}{\bar q(V)}\Big)=\ex(X+Z)=\ex X+\frac{1}{2}
\end{equation}
and
\begin{equation}\label{mom4}
\ex \Big(\frac{V^2}{\bar q(V)}\Big)=\ex(X+Z)^2=\ex X^2 +2 \ex X\ex Z +\ex Z^2=\ex  X^2+ \ex X +\frac{1}{3}.
\end{equation}
Now define the estimators $A_n$ and $B_n$ by
\begin{align}
A_n&=\frac{1}{n}\sum_{i=1}^n\frac{V_i}{q(V_i)}-\frac{1}{2},\\
B_n&=\frac{1}{n}\sum_{i=1}^n\frac{V_i^2}{q(V_i)}-A_n-\frac{1}{3}.
\end{align}
By the equations (\ref{mom3}) and (\ref{mom4}) above we see that $A_n$ and $B_n$ are unbiased estimators of $\ex X$ and $\ex X^2$. Moreover, we can estimate the variance of $X$ by $C_n=B_n-A_n^2$.

For the Beta$(\alpha,\beta)$ distribution we have
\begin{align*}
&\ex X= \frac{\alpha}{\alpha +\beta},\\
&\var X = \frac{\alpha\beta}{(\alpha+\beta)^2(\alpha+\beta+1)}.
\end{align*}
Solving these two equations and plugging in our estimators for the expectation and variance of $X$ we get
\begin{align}
&\hat\alpha=A_n\Big(\frac{A_n(1-A_n)}{C_n}-1\Big),\\
&\hat\beta=(1-A_n)\Big(\frac{A_n(1-A_n)}{C_n}-1\Big).
\end{align}
See for instance  Johnson, Kotz and Balakrishnan (1995) for this method of moments estimation procedure for Beta distributions.

Replacing $F$ and $f$ in the optimal bandwidth (\ref{optband}) by $F_{\hat\alpha,\hat\beta}$ and $f_{\hat\alpha,\hat\beta}$ should give a reasonable bandwidth if the true distribution of $X$ is close to a Beta distribution.

\bigskip

In the situation where the density $q$ is unknown the asymptotic bias is equal to
\begin{equation}\label{bias2}
\frac{1}{2}\,\int_{-\infty}^\infty v^2 w(v)dv\Big(h^2((1-F(x)b^-(x)+F(x)b^+(x))
 -{\tilde h}^2\,\frac{f(x)q''(x)}{q(x)}\Big).
\end{equation}
We can estimate $q$ from the observations times $T_1,\ldots,T_n$, using a bandwidth selector for the optimal bandwidth $\tilde h$ in estimating $q''$ which is of order $n^{-1/7}$. See for instance H\"{a}rdle,  Marron and Wand (1990) for a least squares cross validation method. We can use the resulting estimates of $q, q'$ and $q''$ in the optimal bandwidth (\ref{optband}) where we have to add the extra term in (\ref{bias2}) to the bias. Now we have estimates of
$q, q'$ and $q''$ that we can use in the optimal bandwidth, we have to estimate the distribution of $X$. As above we can use a Beta reference bandwidth. We can estimate the parameters as above but now with estimated   density $q$. Here we use a different bandwidth since only an estimate of $q$ itself is needed. Again a cross validation method can be used, or the Sheather and Jones bandwidth as in Sheather and Jones (1991). these bandwidths will be of order $n^{-1/5}$. Again this method should work fine if the true density of $X$ is close to some Beta density.

\section{A simulated example}\label{simulation}

In this section an example is given to illustrate the estimators found in Section  \ref{const} and \ref{unknown}.
We have simulated the final estimator with the assumption that the density $q$ is known for USCD in which $X_i \sim Beta(2,2)$ and $T_i \sim {\cal N}(0.5,0.3)$  conditioned on $0\leq T_i\leq 1$. The same simulations for the case of unknown density $q$ resulted in graphs of minimal difference compared to the graphs in this section and are therefore not displayed. For more simulated examples we refer to   Graafland (2017).

The kernel we used is equal to the biweight kernel
$$
w(x)=\frac{15}{16}(1-x^2)^2i_{[-1,1]}(x).
$$
The estimator $\hat F_n(x)$ is chosen equal to
$$
\hat F_n(x)=F_n^{1/2}(x)=\frac{1}{2}(F_n^-(x)+F_n^+(x)),
$$
with $h_2 = 0.16$. As shown in Section \ref{chaconsistentF} this estimator satisfies the assumptions of Theorem \ref{theorem2} and \ref{ThmUnkownq}.
Finally, the sample size $n$ equals $10000$ and the bandwidth $h_1$ in $g_{nh}$ equals $0.22$.\\
To reduce computations we have implemented a WARPing technique as described in H\"{a}erdle (1991).

In Figure \ref{fig1} the true density $f$ of the variables $X_i$, the known density $q$ of the variables $T_i$ and the density $g$ of the transformed variables $V_i$ are plotted. In
 Figure \ref{fig2} the left and right estimates $f_{nh}^{-}$ and $f_{nh}^{+}$ are plotted. The graphs of $f_{nh}^{-}$ and $f_{nh}^{+}$ show boundary effects due to the discontinuity of the derivatives of $g$ at $x = 1$.  The final estimate $f_{nh}$ is plotted in Figure \ref{fig3}. Due to the factors $1-\hat{F}_n(x)$ and $\hat{F}_n(x)$ the boundary effects are greatly reduced in the case of our final estimate $f_{nh}$. (The contribution of the estimate $\hat{F}_n(x)$ barely adds new boundary effects as can be seen in Figure \ref{fig4}.) The bias and the variance of the final estimate  $f_{nh}$ are reduced conform Theorem \ref{theorem2} and \ref{ThmUnkownq} that promise asymptotic optimality with respect to the MSE.

In Figure \ref{fig4} the estimates $F_{nh}^{-}$ and $\hat{F}_n =F_{nh}^{1/2}$ are plotted.  Again, the graph of $\hat{F}_n =F_{nh}^{1/2}$ shows a reduction (not necessarily optimal) of the bias and the variance. Boundary effects due to the shape of $g$ occur but are small in both graphs as mentioned before. \\

\begin{rem}{\rm
In our simulated example the density $q$ is nicely bounded away from zero. Simulation of other examples, Graafland (2017), show that the performance is worse when $q$ is near zero. This is of course predicted by the asymptotic variance of our estimator which contains a factor $1/q(x)$.}
\end{rem}
\begin{rem}{\rm
Note that the density $g$ typically has a kink at one. The kernel estimation literature     proposes several methods to correct for boundary effects, see, e.g. Jones (1993). We have not applied these methods in this example as the impact of the boundary effects is small for the final estimator $f_{nh}$.}
\end{rem}
\begin{figure}[h]\label{fig1}
$$
\includegraphics[width=14cm]{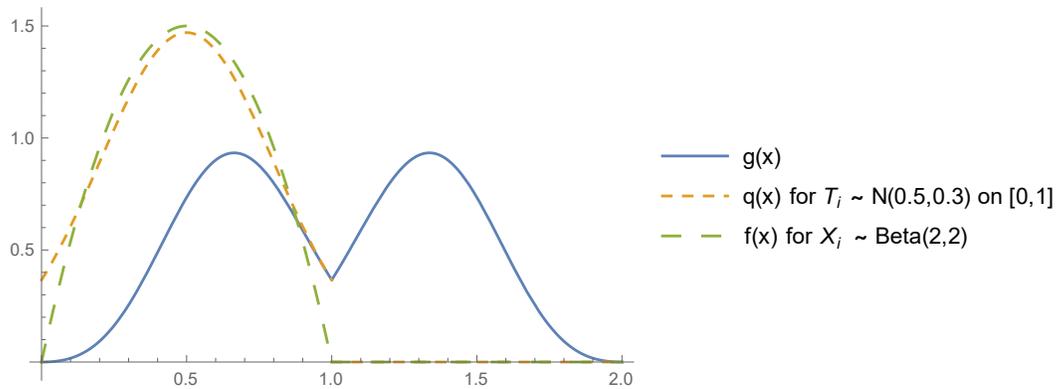}
$$
\caption[]{The density $f$, $q$ and $g$ for  $X_i \sim Beta(2,2)$ and $T_i \sim {\cal N}(0.5,0.3)$  conditioned on $0\leq T_i\leq 1$ .\label{fig1}}
\end{figure}
\begin{figure}[h]\label{fig2}
$$
\includegraphics[width=7.5cm]{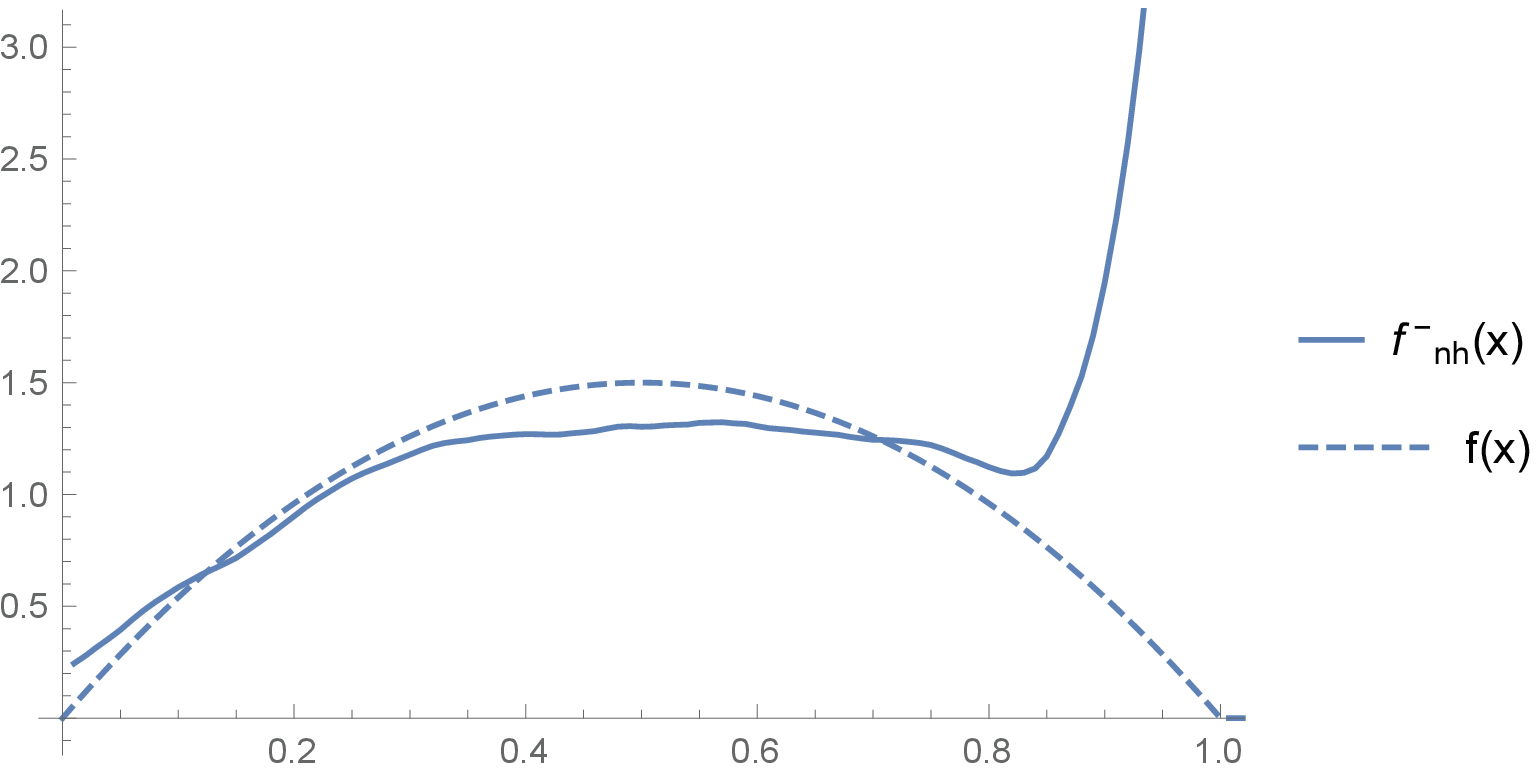}
\includegraphics[width=7.5cm]{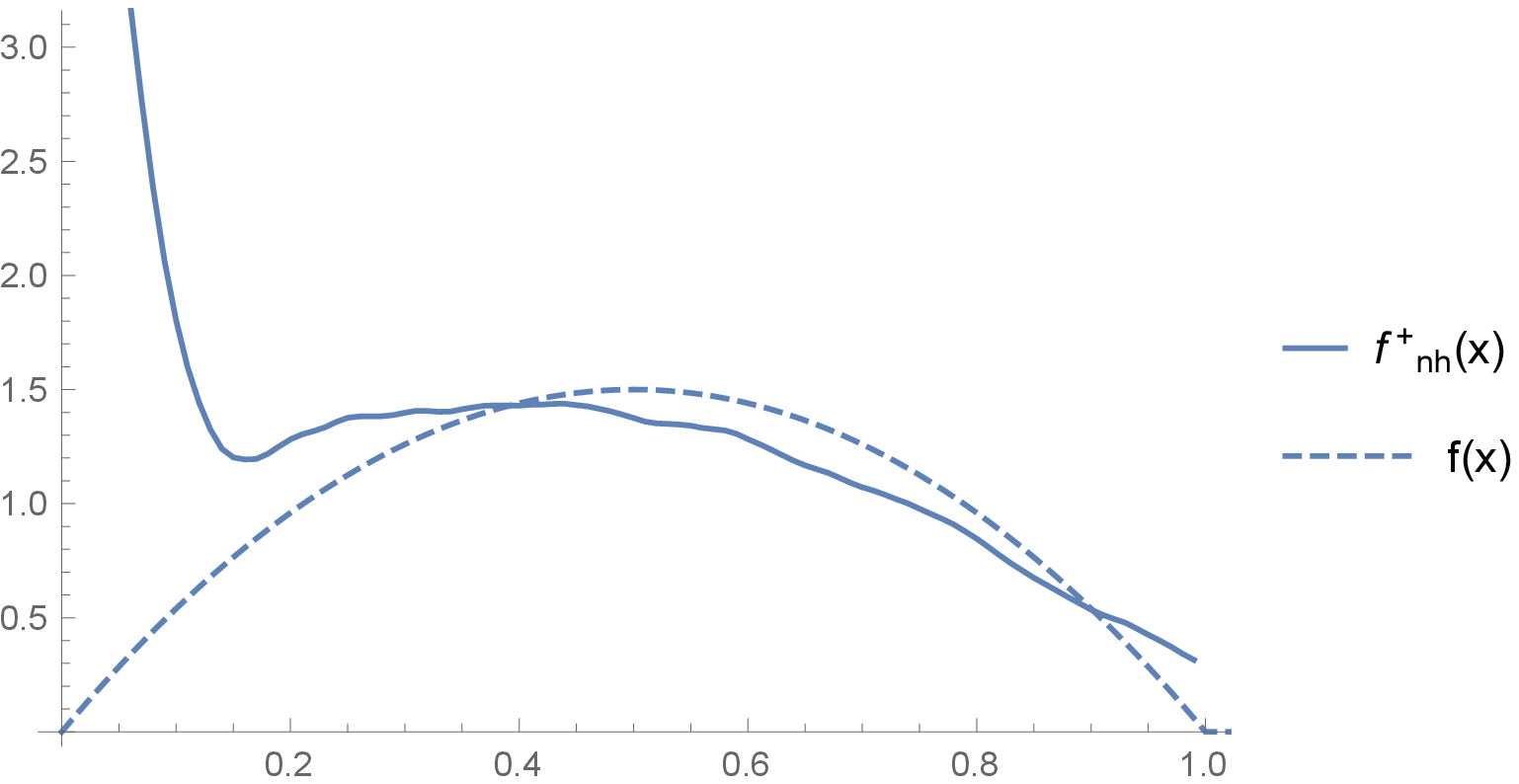}
 $$
\caption[]{Left: the left density estimate $f_{nh}^-$ and the true $f$. Right: the right density estimate $f_{nh}^-$ and the true $f$. \label{fig2} }
\end{figure}
\begin{figure}[h]\label{fig3}
$$
\includegraphics[width=7.5cm]{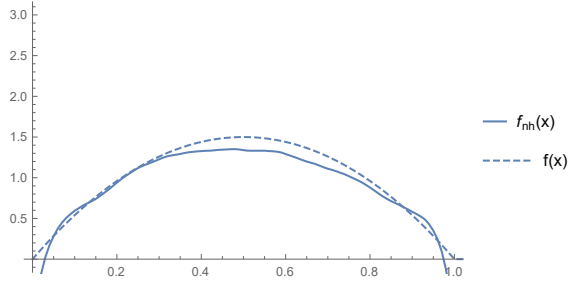}
$$
\caption[]{The final estimate $f_{nh}$ and the true $f$.\label{fig3} }
\end{figure}
\begin{figure}[h]\label{fig4}
$$
\includegraphics[width=7.5cm]{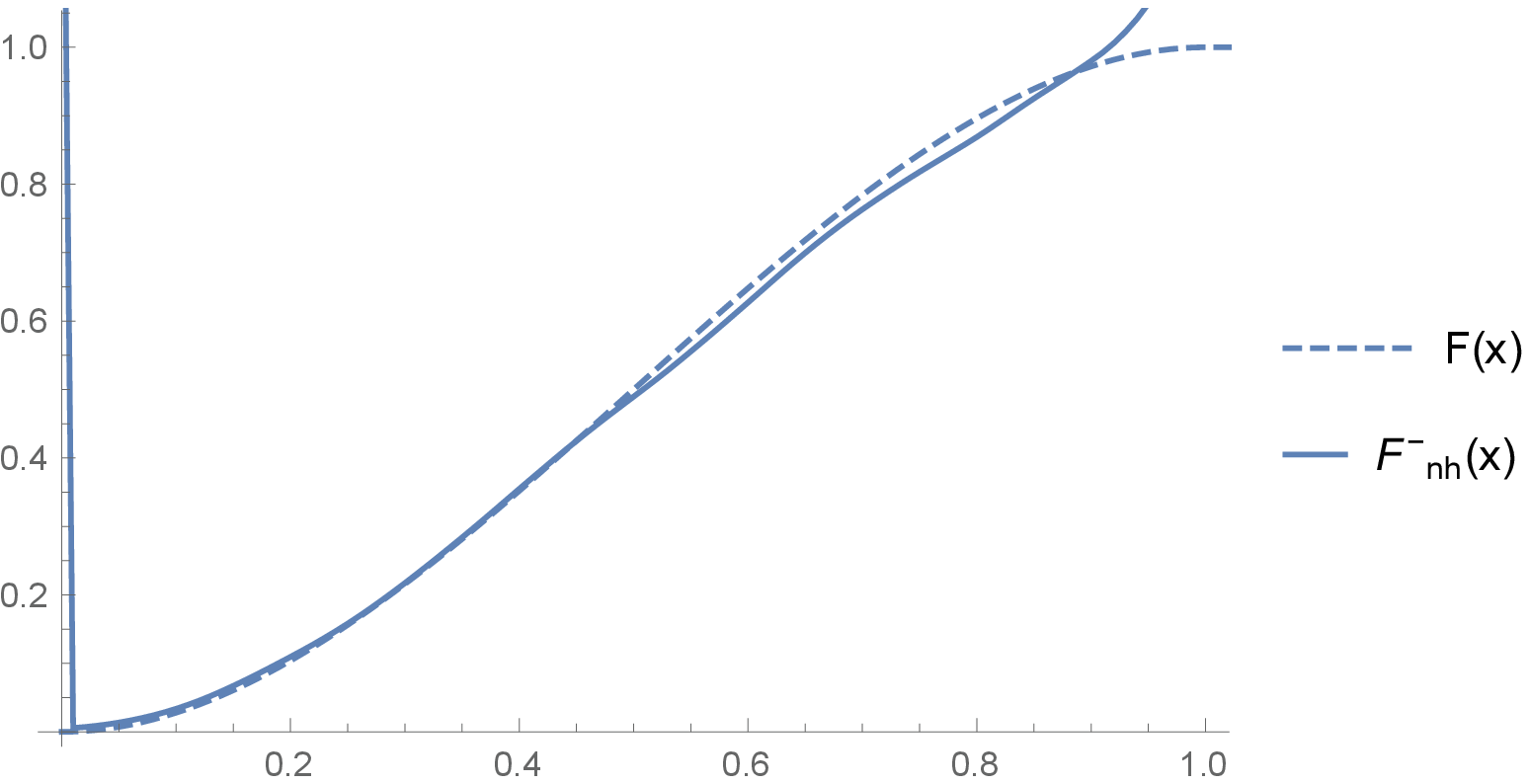}
\includegraphics[width=7.5cm]{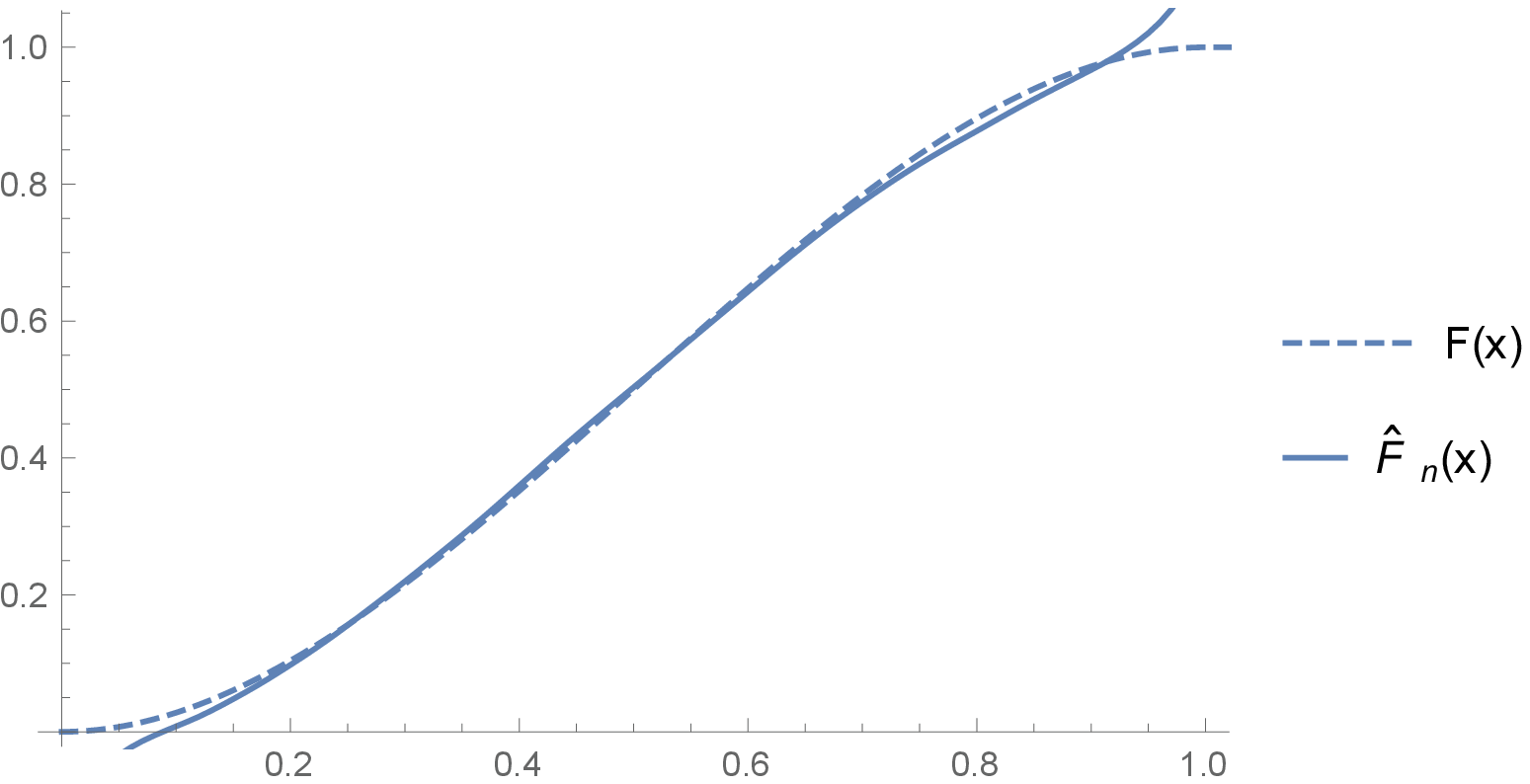}
$$
\caption[]{Left: the left estimate $F_{nh}^-$ and the true $F$. Right: the combined estimate $\hat F_n$ and the true $F$.\label{fig4} }
\end{figure}

\newpage

\section{Proofs}\label{proofs}

\subsection{Proof of Theorem \ref{exvarfnhmin}}
The expansions of the expectation and variance of   $g_{nh}$ and its derivative are standard and given in the following lemma. The proof is omitted.
\begin{lem}\label{ExVargnh}
Assume that the density $g$ is bounded or integrable. Assume also that $g$ is twice differentiable with continuous and bounded $g''$. If the kernel function $w$ satisfies Condition W, then we have
\begin{equation}\label{exgnh}
\begin{split}
\ex g_{nh}(x) = g(x) + \frac{1}{2}h^2g''(x)\int v^2w(v)dv + o(h^2),\\
\var g_{nh}(x) = \frac{1}{nh}g(x)\int w(u)^2du + o\left(\frac{1}{nh}\right).
\end{split}
\end{equation}

If furthermore, $g$ is three times differentiable with continuous and bounded $g^{(3)}$, then we have

\begin{equation}\label{exg'nh}
\begin{split}
\ex g_{nh}'(x) = g'(x) + \frac{1}{2}h^2g^{(3)}(x)\int v^2w(v)dv + o(h^2),\\
\var g_{nh}'(x) = \frac{1}{nh^3}g(x)\int w'(u)^2du + o\left(\frac{1}{nh^3}\right).
\end{split}
\end{equation}
\end{lem}

\noindent{\bf Proof of Theorem \ref{exvarfnhmin}}\\
The assumptions on $f$ and $q$ imply that a function $g$ of the form (\ref{gdens1}) is three times differentiable with continuous and bounded derivative $g^{(3)}$ on the intervals $(0,1)$ and $(1,2)$. Hence the assumptions of Lemma \ref{ExVargnh} are in particular satisfied on the intervals $(0,1)$ for $g$ of the form (\ref{gdens1}). With this in mind we first prove (\ref{exfnhmin}).

 We have
\begin{equation}\nonumber
\begin{split}
\ex f^{-}_{nh}(x) &= \ex \left( \frac{1}{q(x)} g_{nh}'(x) - \frac{q'(x)}{q^{2}(x)}g_{nh}(x)\right) \space  \text{by (\ref{estimator1})}\\
&= \frac{1}{q(x)} \ex g_{nh}'(x) - \frac{q'(x)}{q^{2}(x)}\ex g_{nh}(x)\\
&= \frac{1}{q(x)}\left(g'(x) + \frac{1}{2}h^2g'''(x)\int v^2w(v)dv + o(h^2)\right)\\
& \quad - \frac{q'(x)}{q^{2}(x)}\left(g(x) + \frac{1}{2}h^2g''(x)\int v^2w(v)dv + o(h^2)\right)\\
&= f(x) + \frac{1}{q(x)}\left(\frac{1}{2}h^2g'''(x)\int v^2w(v)dv + o(h^2)\right)\\
& \quad - \frac{q'(x)}{q^{2}(x)}\left(\frac{1}{2}h^2g''(x)\int v^2w(v)dv + o(h^2)\right)\\
&= f(x) + \frac{1}{2}h^2\int v^2w(v)dv\left(\frac{1}{q(x)}g'''(x) - \frac{q'(x)}{q^{2}(x)}g''(x)\right) + o(h^2).
\end{split}
\end{equation}
We continue with the expression for $\var f^{-}_{nh}(x)$. We have
\begin{equation}\nonumber
\begin{split}
\var f^{-}_{nh}(x) &= \var \left( \frac{1}{q(x)} g_{nh}'(x) - \frac{q'(x)}{q^{2}(x)}g_{nh}(x)\right)\\
&= \frac{1}{q(x)^2}\var g_{nh}'(x) - \frac{q'(x)^2}{q^{4}(x)}\var g_{nh}(x) - 2\frac{1}{q(x)}\frac{q'(x)}{q^2(x)}\cov \left(g_{nh}'(x),g_{nh}(x)\right).
\end{split}
\end{equation}
We consider the second and third term in  $\var f^{-}_{nh}(x)$ separately.
Notice that $\var g_{nh}(x) \ll \var g_{nh}'(x)$, because $1/(nh) \ll 1/(nh^3)$. Hence the second term is negligible.
By the Cauchy Schwarz inequality we have
\begin{align*}
\cov (g_{nh}'(x),g_{nh}(x))&\leq \sqrt{\var g_{nh}'(x)}\sqrt{\var g_{nh}(x)}\\
&\ll\sqrt{\var g_{nh}'(x)}\sqrt{\var g_{nh}'(x)}\\
&=\var g'_{nh}(x).
\end{align*}
Hence the third term is negligible.
We may therefore conclude
\begin{equation}\nonumber
\begin{split}
\var f^{-}_{nh}(x) &= \frac{1}{q(x)^2}\var g_{nh}'(x) + o(\var g'_{nh}(x)) \\
&= \frac{1}{q(x)^2}\left(\frac{1}{nh^3}g(x)\int w'(u)^2du + o\left(\frac{1}{nh^3}\right)\right) + o\left(\frac{1}{nh^3}\right)\\
&= \frac{1}{q(x)}\frac{1}{nh^3}F(x)\int w'(u)^2du + o\left(\frac{1}{nh^3}\right).
\end{split}
\end{equation}$\hfill \Box$

\subsection{Proof of Theorem \ref{exvarfnhplus}}
By a similar reasoning as in the proof of Theorem \ref{exvarfnhmin} the assumptions of Lemma \ref{ExVargnh} are now in particular satisfied on the interval $(1,2)$ for $g$ of the form (\ref{gdens1}). With this in mind we first prove equation (\ref{exfnhplus}).

\begin{equation}\nonumber
\begin{split}
\ex f^{+}_{nh}(x) &= \ex \left( -\frac{1}{q(x)} g_{nh}'(x+1) + \frac{q'(x)}{q^{2}(x)}g_{nh}(x+1)\right) \space  \text{by (\ref{estimator2})}\\
&= -\frac{1}{q(x)} \ex g_{nh}'(x+1) + \frac{q'(x)}{q^{2}(x)}\ex g_{nh}(x+1)\\
&= -\frac{1}{q(x)}\left(g'(x+1) + \frac{1}{2}h^2g'''(x+1)\int v^2w(v)dv + o(h^2)\right)\\
& \quad + \frac{q'(x)}{q^{2}(x)}\left(g(x+1) + \frac{1}{2}h^2g''(x+1)\int v^2w(v)dv + o(h^2)\right)\\
&= f(x) + \frac{1}{q(x)}\left(-\frac{1}{2}h^2g'''(x+1)\int v^2w(v)dv + o(h^2)\right)\\
& \quad + \frac{q'(x)}{q^{2}(x)}\left(\frac{1}{2}h^2g''(x+1)\int v^2w(v)dv + o(h^2)\right)\\
&= f(x) + \frac{1}{2}h^2\int v^2w(v)dv\left(-\frac{1}{q(x)}g'''(x+1) + \frac{q'(x)}{q^{2}(x)}g''(x+1)\right) + o(h^2).
\end{split}
\end{equation}
We continue with the expression for $\var f^{+}_{nh}(x)$. We have

\begin{equation}\nonumber
\begin{split}
\var f^{+}_{nh}(x) &= \var \left( -\frac{1}{q(x)} g_{nh}'(x+1) + \frac{q'(x)}{q^{2}(x)}g_{nh}(x+1)\right)\\
&= \frac{1}{q(x)^2}\var g_{nh}'(x+1) + \frac{q'(v)^2}{q^{4}(x)}\var g_{nh}(x+1)\\
& \quad + 2\frac{1}{q(x)}\frac{q'(x)}{q^2(x)}\cov \left(g_{nh}'(x+1),g_{nh}(x+1)\right).
\end{split}
\end{equation}
By the same reasoning as in the proof of Theorem \ref{exvarfnhmin} we have that $\var g_{nh}(x+1) \ll \linebreak \var g_{nh}'(x+1)$. This result and the Cauchy Schwarz inequality give $\cov (g_{nh}'(x+1),g_{nh}(x+1))\ll \linebreak \var g_{nh}'(x+1)$. We may therefore conclude
\begin{equation}\nonumber
\begin{split}
\var f^{+}_{nh}(x) &= \frac{1}{q(x)^2}\var g_{nh}'(x+1) + o(\var g_{nh}'(x+1)) \\
&= \frac{1}{q(x)^2}\left(\frac{1}{nh^3}g(x+1)\int w'(u)^2du + o\left(\frac{1}{nh^3}\right)\right) + o\left(\frac{1}{nh^3}\right)\\
&= \frac{1}{q(x)}\frac{1}{nh^3}(1-F(x))\int w'(u)^2du + o\left(\frac{1}{nh^3}\right).\\
\end{split}
\end{equation}$\hfill\Box$

\subsection{Proof of Theorem \ref{ExVarfnht}}
We first prove (\ref{exfnht}). We apply Theorem \ref{exvarfnhmin} and Theorem \ref{exvarfnhplus}. We have
\begin{equation}\nonumber
\begin{split}
\ex f^{t}_{nh}(x) &= \ex \left( t f^{-}_{nh}(x) + (1-t)f^{+}_{nh}(x) \right)\\
&= t \ex f^{-}_{nh}(x) + (1-t) \ex f^{+}_{nh}(x)\\
&= t\left(f(x) + \frac{1}{2}h^2\int v^2w(v)dv \: b^-(x) + o(h^2)\right)\\
& \quad + (1-t)\left(f(x) + \frac{1}{2}h^2\int v^2w(v)dv \: b^+(x) + o(h^2)\right)\\
&= f(x) + \frac{1}{2}h^2\int v^2w(v)dv \left(tb^-(x) + (1-t)b^+(x)\right) +o(h^2).
\end{split}
\end{equation}
In order to prove (\ref{varfnht}) we use the following lemma to bound the covariance. Its proof is given in Section \ref{lemmaproof}.
\begin{lem}\label{lemmacov} Under the assumptions of Theorem \ref{ExVarfnht} we have
$$\cov \left(f^{-}_{nh}(x),f^{+}_{nh}(x)\right) = o\left(\frac{1}{nh^3}\right).$$
\end{lem}
We continue with the expansion of $\var f^{t}_{nh}(x)$. We have
\begin{equation}\nonumber
\begin{split}
\var f^{t}_{nh}(x) &= \var \left(t f^{-}_{nh}(x) + (1-t)f^{+}_{nh}(x)\right)\\
&= t^2 \var f^{-}_{nh}(x) + (1-t)^2 \var f^{+}_{nh}(x)+ 2t(1-t)\cov \left(f^{-}_{nh}(x),f^{+}_{nh}(x)\right)\\
&= t^2 \var f^{-}_{nh}(x) + (1-t)^2 \var f^{+}_{nh}(x) + 2t(1-t)o\left(\frac{1}{nh^3}\right)\\
&= t^2\left(\frac{1}{q(x)}\frac{1}{nh^3}F(x)\int w'(u)^2du + o\left(\frac{1}{nh^3}\right)\right) \\
& \quad + (1-t)^2\left(\frac{1}{q(x)}\frac{1}{nh^3}(1-F(x))\int w'(u)^2du + o\left(\frac{1}{nh^3}\right)\right)+ o\left(\frac{1}{nh^3}\right)\\
&= \frac{1}{q(x)}\frac{1}{nh^3}\left(t^2F(x) + (1-t)^2(1-F(x)\right)\int w'(u)^2du + o\left(\frac{1}{nh^3}\right),
\end{split}
\end{equation}
which shows (\ref{varfnht}).$\hfill \Box$

\subsection{Proof of Lemma \ref{lemmacov}}\label{lemmaproof}
We have
\begin{equation}\nonumber
\begin{split}
\cov \left(f^{-}_{nh}(x),f^{+}_{nh}(x)\right) = & \\
\cov \Big(\frac{1}{q(x)} g_{nh}'(x)& - \frac{q'(x)}{q^{2}(x)}g_{nh}(x),- \frac{1}{q(x)}g_{nh}'(x+1) + \frac{q'(x)}{q^{2}(x)}g_{nh}(x+1)\Big),
\end{split}
\end{equation}
which we can rewrite as
\begin{align*}
\cov \left(f^{-}_{nh}(x),f^{+}_{nh}(x)\right)&=\cov\left(\frac{1}{q(x)} g_{nh}'(x),- \frac{1}{q(x)}g_{nh}'(x+1)\right)&\,(1)\\
&+\cov\left(\frac{1}{q(x)} g_{nh}'(x),\frac{q'(x)}{q^{2}(x)}g_{nh}(x+1)\right)&\,(2)\\
&+\cov\left(- \frac{q'(x)}{q^{2}(x)}g_{nh}(x),- \frac{1}{q(x)}g_{nh}'(x+1)\right)&\,(3)\\
&+\cov\left(- \frac{q'(x)}{q^{2}(x)}g_{nh}(x),\frac{q'(x)}{q^{2}(x)}g_{nh}(x+1)\right).&\,(4)\\
\end{align*}
We work out the first line separately and show that the second, third and fourth line are negligible compared to the first line. We have

\begin{equation}\nonumber
\begin{split}
(1)\quad \cov\left(\frac{1}{q(x)} g_{nh}'(x),-\frac{1}{q(x)}g_{nh}'(x+1)\right) = &\\
 -\frac{1}{q(x)^2}\big(\ex &g_{nh}'(x)g_{nh}'(x+1) - \ex g_{nh}'(x)\ex g_{nh}'(x+1)\big).
\end{split}
\end{equation}

Now note that, for $n$ large enough,
\begin{align*}
\ex g_{nh}'(x)g_{nh}'(x+1) &= \ex \Bigg(\sum_{i=1}^n \frac{1}{nh^2}\, w'\Big(\frac{x-V_i}{h}\Big)\,\sum_{j=1}^n \frac{1}{nh^2}\, w'\Big(\frac{x+1-V_j}{h}\Big)\Bigg)\\
&=\ex\Bigg(\sum_{i=1}^n \frac{1}{n^2h^4}w'\Big(\frac{x-V_i}{h}\Big)w'\Big(\frac{x+1-V_i}{h}\Big)\\
& \quad + \sum_{i\neq j}^n\frac{1}{n^2h^4}w'\Big(\frac{x-V_i}{h}\Big)w'\Big(\frac{x+1-V_j}{h}\Bigg)\\
&=\ex \sum_{i\neq j}^n\frac{1}{n^2h^4}w'\Big(\frac{x-V_i}{h}\Big)w'\Big(\frac{x+1-V_j}{h}\Big)\\
&=\sum_{i\neq j}^n\frac{1}{n^2h^4}\ex \left(w'\Big(\frac{x-V_i}{h}\Big)w'\Big(\frac{x+1-V_j}{h}\Big)\right)\\
&=\sum_{i\neq j}^n\frac{1}{n^2h^4}\ex w'\Big(\frac{x-V_i}{h}\Big)\ex w'\Big(\frac{x+1-V_j}{h}\Big)\\
&=\frac{n(n-1)}{n^2}\ex g_{nh}'(x)\ex g_{nh}'(x+1)\\
&=(1-\frac{1}{n})\ex g_{nh}'(x)\ex g_{nh}'(x+1)\\
&=\ex g_{nh}'(x)\ex g_{nh}'(x+1) - \frac{1}{n}O(1)\\
&=\ex g_{nh}'(x)\ex g_{nh}'(x+1) - o\left(\frac{1}{nh^3}\right).
\end{align*}
We have used that $w'\Big(\frac{x-V_i}{h}\Big)w'\Big(\frac{x+1-V_i}{h}\Big) = 0$ for all $i$ provided that $h < 1/2$, which is true for $n$ large enough. We have also used that $V_i$ and $V_j$ are independent for all $i,j$ and the fact that $\ex g_{nh}'(x)\ex g_{nh}'(x+1) = g'(x)g'(x+1)+ O(h^2)= O(1)$.\\

We may now conclude
\begin{align*}
\cov&\left(\frac{1}{q(x)} g_{nh}'(x),- \frac{1}{q(x)}g_{nh}'(x+1)\right)\\
&=-\frac{1}{q(x)^2}\left(\ex g_{nh}'(x)\ex g_{nh}'(x+1) - o\left(\frac{1}{nh^3}\right) - \ex g_{nh}'(x)\ex g_{nh}'(x+1)\right)\\
&= o\left(\frac{1}{nh^3}\right)
\end{align*}
and therefore $(1) = o\left(\frac{1}{nh^3}\right)$.\\

With Cauchy Schwarz it is easily seen that $(2)\ll \var g_{nh}'(x)$ and hence $(2) = o(\frac{1}{nh^3})$, $(3)\ll \var g_{nh}'(x+1)$ and hence $(3) = o(\frac{1}{nh^3})$ and the last term $(4)\ll \var  g_{nh}'(x)$ and hence also $(4) = o(\frac{1}{nh^3})$.$\hfill\Box$\\

\bigskip

\subsection{Proof of Theorem \ref{AsNormfnht}}
Write $f^{t}_{nh}(x)$ for $0 < x < 1$ as follows.
\begin{align*}
f^{t}_{nh}(x) &= t f^{-}_{nh}(x) + (1-t)f^{+}_{nh}(x)\\
&= t\left(\frac{1}{q(x)} g_{nh}'(x) - \frac{q'(x)}{q^{2}(x)}g_{nh}(x)\right)\\
& \quad +(1-t) \left(-\frac{1}{q(x)}g_{nh}'(x+1) - \frac{q'(x)}{q^{2}(x)}g_{nh}(x+1) \right)\\
&= t \left(\frac{1}{q(x)} \sum_{i=1}^n \frac{1}{nh^2}\, w'\Big(\frac{x-V_i}{h}\Big) - \frac{q'(x)}{q^{2}(x)}\sum_{i=1}^n \frac{1}{nh}\, w\Big(\frac{x-V_i}{h}\Big)\right)\\
& \quad +(1-t)\left(-\frac{1}{q(x)}\sum_{i=1}^n \frac{1}{nh^2}\, w'\Big(\frac{x+1-V_i}{h}\Big) - \frac{q'(x)}{q^{2}(x)}\sum_{i=1}^n \frac{1}{nh}\, w\Big(\frac{x+1-V_i}{h}\Big)\right)\\
&= \frac{1}{n}\sum_{i=1}^n\,U_{ih}(x)
\end{align*}
with
\begin{align*}
U_{ih} &= t\left(\frac{1}{q(x)} \frac{1}{h^2}\, w'\Big(\frac{x-V_i}{h}\Big) - \frac{q'(x)}{q^{2}(x)} \frac{1}{h}\, w\Big(\frac{x-V_i}{h}\Big)\right)\\
& \quad +(1-t)\left(-\frac{1}{q(x)} \frac{1}{h^2}\, w'\Big(\frac{x+1-V_i}{h}\Big) - \frac{q'(x)}{q^{2}(x)} \frac{1}{h}\, w\Big(\frac{x+1-V_i}{h}\Big)\right).
\end{align*}
Note that
\begin{equation}\label{eq:exftnhexUih}
\ex f^{t}_{nh}(x) = \ex \frac{1}{n}\sum_{i=1}^n\,U_{ih}(x) = \frac{1}{n} \ex \sum_{i=1}^n\,U_{ih}(x) = \frac{1}{n}n \ex \,U_{1h}(x) = \ex U_{1h}(x)
\end{equation}
and
\begin{equation}\label{eq:varftnhvarUih}
\var f^{t}_{nh}(x) = \var \frac{1}{n}\sum_{i=1}^n\,U_{ih}(x) = \frac{1}{n^2} \var \sum_{i=1}^n\,U_{ih}(x) = \frac{1}{n} \var U_{1h}(x).
\end{equation}

We need the following lemma to prove that $f^{t}_{nh}(x) -\ex f^{t}_{nh}(x)$ ($= \frac{1}{n}\sum_{i=1}^n\,\left(U_{1h}(x) - \ex U_{1h}(x)\right)$) is asymptotically normal distributed.
The lemma enables us to show that the Lyapunov condition in the Central Limit Theorem is satisfied. The proof of this lemma can be found in Section \ref{lemmaproof2}. \\

\begin{lem}\label{lem:um}
Under the assumptions of Theorem \ref{AsNormfnht}, we have for $m$ even and for $h \to 0$
\begin{equation}
\ex U_{1h}(x)^m = \frac{1}{q(x)^{m-1}}\frac{1}{h^{2m-1}}\left(t^mF(x) + (1-t)^m(1-F(x)\right)\int w'(u)^mdu + o\left(\frac{1}{h^{2m-1}}\right).
\end{equation}
\end{lem}
\bigskip
We now check the Lyapunov condition with the help of Lemma \ref{lem:um}. This means that for some $\delta >0$ we have to check
$$
\frac{\ex |U_{1h}(x)-\ex
U_{1h}(x)|^{2+\delta}}{n^{\delta/2}(\var(U_{1h}(x)))^{1+\delta/2}}\to 0.
$$

We use that $\ex U_{1h}(x) = \ex f^{t}_{nh}(x) = O(1)$ and $\var U_{1h}(x) =  O(1/nh^3)$. Furthermore we use that $(a+b)^4\leq 2^3(a^4+b^4)$. For $\delta =2$ we have

\begin{align*}
\frac{\ex |U_{1h}(x)-\ex
U_{1h}(x)|^{4}}{n(\var(U_{1h}(x)))^2}&
\leq \frac{2^3(\ex
U_{1h}(x)^4+(\ex U_{1h}(x))^4)}{n(\var(U_{1h}(x)))^2}\\
&\sim \frac{8\left(\frac{1}{q(x)^{3}}\frac{1}{h^{7}}\left(t^4F(x) + (1-t)^4(1-F(x)\right)\int w'(u)^4du + o(\frac{1}{h^{7}})\right)}{n(\frac{1}{h^3}c_2)^2}\\
&\sim \frac{8c_1}{nhc_2^2}\to 0,
\end{align*}
as $n \to 0.$ Hence the Lyapunov condition is satisfied which proves the theorem. $\hfill\Box$

\subsection{ Proof of Lemma \ref{lem:um}}\label{lemmaproof2}
We have
\begin{align*}
\ex U_{1h}(x)^m &= \ex \Bigg( t\left(\frac{1}{q(x)} \frac{1}{h^2}\, w'\Big(\frac{x-V_i}{h}\Big) - \frac{q'(x)}{q^{2}(x)} \frac{1}{h}\, w\Big(\frac{x-V_i}{h}\Big)\right)\\
& \quad +(1-t)\left(-\frac{1}{q(x)} \frac{1}{h^2}\, w'\Big(\frac{x+1-V_i}{h}\Big) - \frac{q'(x)}{q^{2}(x)} \frac{1}{h}\, w\Big(\frac{x+1-V_i}{h}\Big) \right)\Bigg)^m.
\end{align*}
Write $\ex U_{1h}(x)^m = \ex [(a+b)^m]$ with $a = t\left(\frac{1}{q(x)} \frac{1}{h^2}\, w'\Big(\frac{x-V_i}{h}\Big) - \frac{q'(x)}{q^{2}(x)} \frac{1}{h}\, w\Big(\frac{x-V_i}{h}\Big)\right)$\\
and $b=(1-t)\left(-\frac{1}{q(x)} \frac{1}{h^2}\, w'\Big(\frac{x+1-V_i}{h}\Big) - \frac{q'(x)}{q^{2}(x)} \frac{1}{h}\, w\Big(\frac{x+1-V_i}{h}\Big) \right).$
Now note for $h < 1/2$ we have
\begin{align*}
w'\Big(\frac{x-V_i}{h}\Big)w'\Big(\frac{x+1-V_i}{h}\Big) &= 0,\\
w'\Big(\frac{x-V_i}{h}\Big)w\Big(\frac{x+1-V_i}{h}\Big) &= 0,\\
w\Big(\frac{x-V_i}{h}\Big)w'\Big(\frac{x+1-V_i}{h}\Big) &= 0,\\
w\Big(\frac{x-V_i}{h}\Big)w\Big(\frac{x+1-V_i}{h}\Big) &= 0.\\
\end{align*}
Hence we have
\begin{align*}
\ex U_{1h}(x)^m &= \ex [(a+b)^m]\\
&= \ex \Big[ \sum_{k=0}^m {{m}\choose{k}} a^{m-k}b^{k} \Big]\\
&= \ex \Big[ a^m + b^m \Big] \qquad  \mbox{because} \qquad a^{m-k}b^k = 0\; \forall \; k\in \{1,\dots, m-1\}\\
&= \ex \Big[ t^m\left(\frac{1}{q(x)} \frac{1}{h^2}\, w'\Big(\frac{x-V_i}{h}\Big) - \frac{q'(x)}{q^{2}(x)} \frac{1}{h}\, w\Big(\frac{x-V_i}{h}\Big)\right)^m\\
& \quad + (1-t)^m(-1)^m\left(\frac{1}{q(x)} \frac{1}{h^2}\, w'\Big(\frac{x+1-V_i}{h}\Big) - \frac{q'(x)}{q^{2}(x)} \frac{1}{h}\, w\Big(\frac{x+1-V_i}{h}\Big) \right)^m\Big].\\
\end{align*}
The fact that
\begin{align*}
\ex \Big[ t^m\Big(\frac{1}{q(x)} \frac{1}{h^2}\, &w'\Big(\frac{x-V_i}{h}\Big) - \frac{q'(x)}{q^{2}(x)} \frac{1}{h}\, w\Big(\frac{x-V_i}{h}\Big)\Big)^m\Big]\\
&\sim t^m \frac{1}{q(x)^{m-1}}\frac{1}{h^{2m-1}}F(x)\int w'(u)^mdu + o\left(\frac{1}{h^{2m-1}}\right)\qquad (1)\\
\end{align*}
and
\begin{align*}
\ex &\Big[ (1-t)^m(-1)^m\Big(\frac{1}{q(x)} \frac{1}{h^2}\, w'\Big(\frac{x+1-V_i}{h}\Big) - \frac{q'(x)}{q^{2}(x)} \frac{1}{h}\, w\Big(\frac{x+1-V_i}{h}\Big) \Big)^m\Big]\\
&\sim (1-t)^m(-1)^m \frac{1}{q(x)^{m-1}}\frac{1}{h^{2m-1}}(1-F(x))\int w'(u)^mdu + o\left(\frac{1}{h^{2m-1}}\right)\qquad (2)
\end{align*}
leads us to the desired result of Lemma \ref{lem:um}.\\

What rests is proving $(1)$ and $(2)$. For $(1)$ we have
\begin{align*}
\ex \Big[ t^m\Big(\frac{1}{q(x)} \frac{1}{h^2}\, &w'\Big(\frac{x-V_i}{h}\Big) - \frac{q'(x)}{q^{2}(x)} \frac{1}{h}\, w\Big(\frac{x-V_i}{h}\Big)\Big)^m\Big]\\
&= t^m \ex \Big[ \sum_{k=0}^m {{m}\choose{k}} \Big(\frac{1}{q(x)}\frac{1}{h^2}\, w'\Big(\frac{x-V_i}{h}\Big)\Big)^{m-k}\Big(- \frac{q'(x)}{q^{2}(x)} \frac{1}{h}\, w\Big(\frac{x-V_i}{h}\Big)\Big)^{k}\Big]\\
&= t^m \sum_{k=0}^m {{m}\choose{k}} \ex \Big[\Big(\frac{1}{q(x)}\frac{1}{h^2}\, w'\Big(\frac{x-V_i}{h}\Big)\Big)^{m-k}\Big(- \frac{q'(x)}{q^{2}(x)} \frac{1}{h}\, w\Big(\frac{x-V_i}{h}\Big)\Big)^{k}\Big].\\
\end{align*}
We next expand  the terms of the summation. We have
\begin{align*}
\ex \Big[\Big(\frac{1}{q(x)}\frac{1}{h^2}\, &w'\Big(\frac{x-V_i}{h}\Big)\Big)^{m-k}\Big(- \frac{q'(x)}{q^{2}(x)} \frac{1}{h}\, w\Big(\frac{x-V_i}{h}\Big)\Big)^{k}\Big]\\
&= \frac{1}{q(x)^{m-k}}\frac{1}{h^{2m-2k}}\Big(- \frac{q'(x)}{q^{2}(x)}\Big)^{k}\frac{1}{h^k} \ex \Big[w'\Big(\frac{x-V_i}{h}\Big)^{m-k}w\Big(\frac{x-V_i}{h}\Big)^{k}\Big]\\
&=\frac{1}{q(x)^{m-k}}\frac{1}{h^{2m-k}}\Big(- \frac{q'(x)}{q^{2}(x)}\Big)^{k} \int w'\Big(\frac{x-t}{h}\Big)^{m-k}w\Big(\frac{x-t}{h}\Big)^{k}g(t)dt\\
&=\frac{1}{q(x)^{m-k}}\frac{1}{h^{2m-k-1}}\Big(- \frac{q'(x)}{q^{2}(x)}\Big)^{k} \int w'(-v)^{m-k}w(-v)^{k}g(x+hv)dv\\
&=\frac{1}{q(x)^{m-k}}\frac{1}{h^{2m-k-1}}\Big(- \frac{q'(x)}{q^{2}(x)}\Big)^{k} \left(\int w'(-v)^{m-k}w(-v)^{k}g(x)dv + o(1)\right)\\
&=\frac{1}{q(x)^{m-k}}\frac{1}{h^{2m-k-1}}\Big(- \frac{q'(x)}{q^{2}(x)}\Big)^{k} g(x)\int w'(-v)^{m-k}w(-v)^{k} dv + o\left(\frac{1}{h^{2m-k-1}}\right)\\
&=(-1)^k\frac{1}{q(x)^{m+k-1}}\frac{1}{h^{2m-k-1}}q'(x)^k F(x)\int w'(-v)^{m-k}w(-v)^{k} dv + o\left(\frac{1}{h^{2m-k-1}}\right).
\end{align*}
In the fourth equation we used the dominated convergence theorem.
We conclude that terms with $k \geq 1$ in the summation   all have order smaller or equal to $o(\frac{1}{h^{2m-1}})$. Hence we may conclude for (1),

\begin{align*}
\ex \Big[ t^m\Big(\frac{1}{q(x)} \frac{1}{h^2}\, &w'\Big(\frac{x-V_i}{h}\Big) - \frac{q'(x)}{q^{2}(x)} \frac{1}{h}\, w\Big(\frac{x-V_i}{h}\Big)\Big)^m\Big]\\
&\sim t^m \frac{1}{q(x)^{m-1}}\frac{1}{h^{2m-1}}F(x)\int w'(u)^mdu + o\left(\frac{1}{h^{2m-1}}\right).
\end{align*}

Statement $(2)$ can be proven in a similar way. $\hfill\Box$

\subsection{Proof of Theorem \ref{theorem2}}
Note that Theorem \ref{AsNormfnht} proves asymptotic normality for the estimator $f_{nh}^t(x)$ with $t = 1-F(x)$, because we have $(1-F(x)) \in [0,1]$. We use this fact to prove asymptotic normality for our final estimator $f_{nh}(x) = f^{1-\hat{F}_{n}(x)}_{nh}(x)$. \\

Write
\begin{equation}\label{decomp}
f_{nh}(x)=(1-\hat F_n(x)) f^{-}_{nh}(x) +  \hat
F_n(x)f^{+}_{nh}(x) =f_{nh}^{1-F(x)}(x) +R_{nh}(x),
\end{equation}
where
\begin{equation}\label{RandWdefs}
R_{nh}(x)=(\hat F_n(x)-F(x)) S_{nh}(x)\quad\mbox{and}\quad
S_{nh}(x)=f^{+}_{nh}(x) -f^{-}_{nh}(x).
\end{equation}
Note that
\begin{align*}
\sqrt{nh^3}(f_{nh}(x)-\ex f_{nh}(x)) &= \sqrt{nh^3} \left(f_{nh}^{1-F(x)}(x) - \ex f_{nh}^{1-F(x)}(x)+ R_{nh}(x)- \ex R_{nh}(x)\right)\nonumber\\
&=\sqrt{nh^3} \left(f_{nh}^{1-F(x)}(x) - \ex f_{nh}^{1-F(x)}(x)\right)\\
&\quad + \sqrt{nh^3}R_{nh}(x) - \sqrt{nh^3}\ex R_{nh}(x).
\end{align*}
We show in the sequel that $\sqrt{nh^3}R_{nh}(x)$ and $\sqrt{nh^3}\ex R_{nh}(x)$ converge to $0$ in distribution.\\

We first rewrite $S_{nh}(x)$, we have
\begin{align*}
S_{nh}(x) &= - \left( \frac{1}{q(x)}g_{nh}'(x+1) - \frac{q'(x)}{q^{2}(x)}g_{nh}(x+1) \right) - \left( \frac{1}{q(x)}g_{nh}'(x) - \frac{q'(x)}{q^{2}(x)}g_{nh}(x) \right)\\
&= - \left( \frac{1}{q(x)}\frac{1}{n}\,\sum_{i=1}^n \frac{1}{h^2}\, w'\Big(\frac{x+1-V_i}{h}\Big) - \frac{q'(x)}{q^{2}(x)}\frac{1}{n}\,\sum_{i=1}^n \frac{1}{h}\, w\Big(\frac{x+1-V_i}{h}\Big)\right)\\
& \quad -\left( \frac{1}{q(x)}\frac{1}{n}\,\sum_{i=1}^n \frac{1}{h^2}\, w'\Big(\frac{x-V_i}{h}\Big)- \frac{q'(x)}{q^{2}(x)}\frac{1}{n}\,\sum_{i=1}^n \frac{1}{h}\, w\Big(\frac{x-V_i}{h}\Big)\right)\\
&=\frac{1}{n}\,\sum_{i=1}^n \Big(- \left( \frac{1}{q(x)}\frac{1}{h^2}\, w'\Big(\frac{x+1-V_i}{h}\Big) - \frac{q'(x)}{q^{2}(x)}\frac{1}{h}\, w\Big(\frac{x+1-V_i}{h}\Big)\right)\\
& \quad -\left( \frac{1}{q(x)}\frac{1}{h^2}\, w'\Big(\frac{x-V_i}{h}\Big)- \frac{q'(x)}{q^{2}(x)} \frac{1}{h}\, w\Big(\frac{x-V_i}{h}\Big)\right)\Big)\\
&= \frac{1}{n}\,\sum_{i=1}^n W_{ih}(x)
\end{align*}
where
\begin{align*}
W_{ih}(x)= \Big(- \left( \frac{1}{q(x)}\frac{1}{h^2}\, w'\Big(\frac{x+1-V_i}{h}\Big) - \frac{q'(x)}{q^{2}(x)}\frac{1}{h}\, w\Big(\frac{x+1-V_i}{h}\Big)\right)\\
-\left( \frac{1}{q(x)}\frac{1}{h^2}\, w'\Big(\frac{x-V_i}{h}\Big)- \frac{q'(x)}{q^{2}(x)} \frac{1}{h}\, w\Big(\frac{x-V_i}{h}\Big)\right)\Big).\\
\end{align*}

The next lemma establishes some properties of
$S_{nh}(x)$. The proof of Lemma \ref{Wlemma} can be found in Section \ref{chaWlemma}.

\begin{lem}\label{Wlemma}
Under de assumptions of Theorem \ref{theorem2} we have
\begin{equation}
\ex S_{nh}(x)=\frac{1}{2}h^2\int v^2w(v)dv(b^+(x) - b^-(x)) + o(h^2),
\end{equation}
\begin{equation}
\ex W_{ih}(x)^m=\frac{1}{q(x)^{m-1}}\frac{1}{h^{2m-1}}\int w'(u)^mdu + o\left(\frac{1}{h^{2m-1}}\right),
\end{equation}
and
\begin{equation}
\sqrt{nh^3}\left(S_{nh}(x) - \ex S_{nh}(x)\right)\convd
N\Big(0,\frac{1}{q(x)}\int w'(u)^2du\Big).\\
\end{equation}
\end{lem}

\bigskip

We start to analyse the term $\sqrt{nh^3}R_{nh}(x)$ and we prove that it converges to zero in distribution. We rewrite the term as
\begin{align*}
\sqrt{nh^3}R_{nh}(x) &= \sqrt{nh^3}(\hat F_n(x)-F(x)) S_{nh}(x)\\
&= \sqrt{nh^3}(\hat F_n(x)-F(x)) (S_{nh}(x) - \ex  S_{nh}(x))\\
&\quad + \sqrt{nh^3}(\hat F_n(x)-F(x))\ex  S_{nh}(x).
\end{align*}

We estimate the first and second term in the last line separately.
By condition (\ref{fhatcond}) we have that  $\hat F_n(x)-F(x)\convp 0$ and hence by Slutsky's Theorem and the result of Lemma \ref{Wlemma} we may conclude \\
$$\sqrt{nh^3}\left(\hat F_n(x)-F(x)\right)\left(S_{nh}(x) - \ex  S_{nh}(x)\right)\convd 0.$$
Furthermore we have by Lemma \ref{Wlemma} that $\ex  S_{nh}(x)$ equals $ O(h^2)$.
Hence for $h = O(n^{-1/7}) $ we have

\begin{equation}\label{ExpRest}
\begin{split}
\sqrt{nh^3}\ex S_{nh}(x) &= O(n^{1/2}h^{7/2})\\
& = O(n^{1/2}(n^{-1/7})^{7/2})\\
& = O(1).\\
\end{split}
\end{equation}
Using again condition (\ref{fhatcond}) we conclude
\begin{equation}\label{ExpRest2}
\sqrt{nh^3}(\hat F_n(x)-F(x))\ex  S_{nh}(x) \convp 0.
\end{equation}
Together equations (\ref{ExpRest}) and (\ref{ExpRest2}) ensure that
$$\sqrt{nh^3}R_{nh}(x) \convd 0.$$\\

We now analyse the term $\sqrt{nh^3}\ex R_{nh}(x)$ and we prove that it converges to zero in distribution as well.
By the Cauchy-Schwarz
inequality we have
\begin{equation}\nonumber
\begin{split}
\ex \sqrt{nh^3}|R_{nh}(x)| &\leq \sqrt{nh^3} (\ex (\hat
F_n(x)-F(x))^2)^{1/2} (\ex (S_{nh}(x))^2)^{1/2}\\
& = \ex (\hat F_n(x)-F(x))^2)^{1/2}O(1)\to 0.
\end{split}
\end{equation}

By the fact that $\sqrt{nh^3}R_{nh}(x)$ and $\sqrt{nh^3}\ex R_{nh}(x)$ both converge  to zero in distribution, we may conclude that $\sqrt{nh^3}(f_{nh}(x)-\ex f_{nh}(x))$ has the
same asymptotic normal distribution as
$\sqrt{nh^3}(f_{nh}^{1-F(x)}(x)-\ex f_{nh}^{1-F(x)}(x))$,
which proves the first statement of the theorem.\\

\bigskip

The second statement of the theorem is proven as follows. By equation (\ref{exfnht}) we have
\begin{equation}\label{exfnhFx}
\ex f_{nh}^{1-F(x)}(x)
= f(x) + \frac{1}{2}h^2\int v^2w(v)dv\left((1-F(x))b^-(x) + F(x)b^+(x)\right)+ o(h^2).
\end{equation}
Furthermore we have
\begin{align}
\ex |R_{nh}(x)| &\leq  (\ex (\hat F_n(x)-F(x))^2)^{1/2}
(\ex (S_{nh}(x))^2)^{1/2} \nonumber\\
&= o\Big(\sqrt{nh^7}\Big)O\Big(\frac{1}{\sqrt{nh^3}}\Big)=o(h^2).\label{biasest}
\end{align}
Together equation (\ref{exfnhFx}) and
(\ref{biasest}) prove the second statement of the theorem. \\

\bigskip

Finally, by equation (\ref{varfnht}), we find that
\begin{align*}\label{varfnhFx}
\var f^{1-F(x)}_{nh}(x) &= \frac{1}{q(x)}\frac{1}{nh^3}\left((1-F(x))^2F(x) + F(x)^2(1-F(x)\right)\int w'(u)^2du + o\left(\frac{1}{nh^3}\right)\\
&=\frac{1}{q(x)}\frac{1}{nh^3}F(x)(1-F(x))\int w'(u)^2du + o\left(\frac{1}{nh^3}\right).
\end{align*}
Thus for sub optimal $h$ we have $\var f_{nh}(x) \sim \var f_{nh}^{1-F(x)}(x)$.
Hence for sub optimal $h$ the estimator $f_{nh}(x)$ is t-optimal.\\
The last line follows by the fact that for $nh^3 \to \infty$ and $h \ll n^{-1/7}$ we have
$$
\sqrt{nh^3}\left(\ex f_{nh}(x) - f(x)\right) \to 0
$$
as $n \to \infty$. $\hfill \Box$\\

Note that Lemma \ref{Wlemma} reveals that the expectation of the difference between the left and right estimator, i.e. $\ex S_{nh}(x)$, depends only on the density $q$ as we have $$
b^+(x) - b^-(x) = \frac{1}{q(x)}q'''(x) - \frac{q'(x)}{q^{2}(x)}q''(x).
$$
This follows from the relation $g(x) + g(x+1)  =  q(x)$ which is evident from the inversion formulas (\ref{inversion1}) and (\ref{inversion2}). Taking derivatives gives similar relations for the higher derivatives.

\begin{rem}{\rm
Note that, in contrast to Theorem \ref{AsNormfnht}, we already assume sufficient smoothness on the functions $f$ and $q$ to prove the first statement of asymptotic normality. In equation (\ref{ExpRest}) this assumption ensures that the difference between $f^+_{nh}(x)$ and $f^-_{nh}(x)$, denoted by $S_{nh}(x)$, satisfies $\ex S_{nh}(x) = O(h^2)$. Together with the restriction $h=O(n^{-1/7})$ on $h$ we obtain $\sqrt{nh^3}\ex S_{nh}(x) = O(1)$.}
\end{rem}

\subsection{Proof of Lemma \ref{Wlemma}}\label{chaWlemma}

The first statement follows from (\ref{exfnhmin}) and (\ref{exfnhplus}).\\
The expresion for $\ex W_{ih}(x)^m$ is obtained as follows,
\begin{align}
\ex W_{ih}(x)^m &= \ex \Bigg(\left(-\frac{1}{q(x)} \frac{1}{h^2}\, w'\Big(\frac{x+1-V_i}{h}\Big) + \frac{q'(x)}{q^{2}(x)} \frac{1}{h}\, w\Big(\frac{x+1-V_i}{h}\Big) \right)\\
&- \left(\frac{1}{q(x)} \frac{1}{h^2}\, w'\Big(\frac{x-V_i}{h}\Big) - \frac{q'(x)}{q^{2}(x)} \frac{1}{h}\, w\Big(\frac{x-V_i}{h}\Big)\right)\Bigg)^m.
\end{align}
Write $\ex W_{ih}(x)^m = \ex [(a+b)^m]$ with $a = -\frac{1}{q(x)} \frac{1}{h^2}\, w'\Big(\frac{x+1-V_i}{h}\Big) + \frac{q'(x)}{q^{2}(x)} \frac{1}{h}\, w\Big(\frac{x+1-V_i}{h}\Big)$\\
and $b=\frac{1}{q(x)} \frac{1}{h^2}\, w'\Big(\frac{x-V_i}{h}\Big) - \frac{q'(x)}{q^{2}(x)} \frac{1}{h}\, w\Big(\frac{x-V_i}{h}\Big).$
Now note if $h < 1/2$, we have
\begin{align*}
w'\Big(\frac{x-V_i}{h}\Big)w'\Big(\frac{x+1-V_i}{h}\Big) &= 0,\\
w'\Big(\frac{x-V_i}{h}\Big)w\Big(\frac{x+1-V_i}{h}\Big) &= 0,\\
w\Big(\frac{x-V_i}{h}\Big)w'\Big(\frac{x+1-V_i}{h}\Big) &= 0,\\
w\Big(\frac{x-V_i}{h}\Big)w\Big(\frac{x+1-V_i}{h}\Big) &= 0.\\
\end{align*}
Hence we have
\begin{align*}
\ex W_{ih}(x)^m &= \ex [(a+b)^m]\\
&= \ex \Big[ \sum_{k=0}^m {{m}\choose{k}} a^{m-k}b^{k} \Big]\\
&= \ex \Big[ a^m + b^m \Big] \qquad  \mbox{because} \qquad a^{m-k}b^k = 0\; \forall \, k\in [1,m-1]\\
&= \ex \Big[\left(-\frac{1}{q(x)} \frac{1}{h^2}\, w'\Big(\frac{x+1-V_i}{h}\Big) + \frac{q'(x)}{q^{2}(x)} \frac{1}{h}\, w\Big(\frac{x+1-V_i}{h}\Big) \right)^m\\
& \quad + (-1)^m\left(\frac{1}{q(x)} \frac{1}{h^2}\, w'\Big(\frac{x-V_i}{h}\Big) - \frac{q'(x)}{q^{2}(x)} \frac{1}{h}\, w\Big(\frac{x-V_i}{h}\Big)\right)^m\Big].
\end{align*}
The fact that
\begin{align*}
\ex \Big[\Big(-\frac{1}{q(x)} \frac{1}{h^2}\, &w'\Big(\frac{x+1-V_i}{h}\Big) + \frac{q'(x)}{q^{2}(x)} \frac{1}{h}\, w\Big(\frac{x+1-V_i}{h}\Big) \Big)^m\Big]\\
&\sim (-1)^m \frac{1}{q(x)^{m-1}}\frac{1}{h^{2m-1}} F(x)\int w'(u)^mdu + o\left(\frac{1}{h^{2m-1}}\right)\qquad (1)\\
\end{align*}
and
\begin{align*}
\ex \Big[ (-1)^m\Big(\frac{1}{q(x)} \frac{1}{h^2}\, &w'\Big(\frac{x-V_i}{h}\Big) - \frac{q'(x)}{q^{2}(x)} \frac{1}{h}\, w\Big(\frac{x-V_i}{h}\Big)\Big)^m\Big]\\
&\sim (-1)^m \frac{1}{q(x)^{m-1}}\frac{1}{h^{2m-1}}(1-F(x))\int w'(u)^mdu + o\left(\frac{1}{h^{2m-1}}\right)\qquad (2)
\end{align*}
can be proven in the same way (leave out the $t$-depended terms) as fact (1) and (2) in the proof of Lemma \ref{lem:um}.\\
Together (1) en (2) lead to the result of the second statement,
\begin{equation}\nonumber
\ex W_{ih}(x)^m=\frac{1}{q(x)^{m-1}}\frac{1}{h^{2m-1}}(-1)^m\int w'(u)^mdu + o\left(\frac{1}{h^{2m-1}}\right).
\end{equation}

For the third statement we prove that the Lyapunov condition holds.
Note that
$\ex W_{ih}(x) = \ex S_{nh}(x) = O(1)$ and $\var W_{ih}(x) = n\var S_{nh}(x)$. By the second statement we have
$$
\var W_{ih}(x) \sim \frac{1}{q(x)}\frac{1}{h^3}\int w'(u)^2du + o\left(\frac{1}{h^3}\right).
$$
Further use that $(a+b)^4\leq 2^3(a^4+b^4)$. We check the Lyapunov condition for $\delta =2$. We have
\begin{align*}
\dfrac{\ex |W_{ih}(x)-\ex
W_{ih}(x)|^{4}}{n(\var W_{ih}(x)))^2}
&\leq \dfrac{2^3(\ex
W_{ih}(x)^4+(\ex W_{ih}(x))^4)}{n(\frac{1}{nh^3}c_2)^2}\\
&\sim \dfrac{8\left(\frac{1}{q(x)^4}\frac{1}{h^{7}}\int w'(u)^4du + O(1)^4\right)}{n(\frac{1}{h^3}c_2)^2}
\sim \dfrac{c_1}{nhc_2^2}\to 0.
\end{align*}
$\hfill\Box$

\subsection{Proof of Theorem \ref{Theorem3}}\label{chaTheorem3}

The assumptions of Theorem \ref{theorem2} satisfy in particular the assumptions of Lemma \ref{ExVargnh}. Hence we may use Lemma \ref{ExVargnh} to compute expansions for the expectations of $F_{nh}^{-}(x)$ and $F_{nh}^{+}(x)$. We have for $0<x<1$,
\begin{align*}
\ex F_{nh}^{-}(x) &= \frac{1}{q(x)} \ex g_{nh}(x) \\
&= \frac{1}{q(x)}\left(g(x) + \frac{1}{2}h^2g''(x)\int v^2w(v)dv + o(h^2)\right)\\
&= F(x) + \frac{1}{q(x)}\frac{1}{2}h^2g''(x)\int v^2w(v)dv + o(h^2)
\end{align*}
and
\begin{align*}
\ex F_{nh}^{+}(x) &= 1 - \frac{1}{q(x)} \ex g_{nh}(x+1) \\
&= 1 - \frac{1}{q(x)}\left(g(x+1) + \frac{1}{2}h^2g''(x+1)\int v^2w(v)dv + o(h^2)\right)\\
&= F(x) - \frac{1}{q(x)}\frac{1}{2}h^2g''(x+1)\int v^2w(v)dv - o(h^2).
\end{align*}
For $0<x<1$ the expansion of the expectation of $F_{nh}^{t}(x)$ in equation (\ref{a8}) is now obtained as follows.
\begin{equation}
\begin{split}
\ex F_{nh}^{t}(x) &= t\ex F_{nh}^{-}(x) + (1-t)\ex F_{nh}^{+}(x)\\
&= t\left(F(x) + \frac{1}{q(x)}\frac{1}{2}h^2g''(x)\int v^2w(v)dv + o(h^2)\right)\\
&\quad+(1-t)\left(F(x) - \frac{1}{q(x)}\frac{1}{2}h^2g''(x+1)\int v^2w(v)dv - o(h^2)\right)\\
&= F(x) + \frac{1}{q(x)}\frac{1}{2}h^2\left(tg''(x)-(1-t)g''(x+1)\right)\int v^2w(v)dv + o(h^2).
\end{split}
\end{equation}
To obtain the expansion for the variance as stated in equation (\ref{a9}) we rewrite $F^{t}_{nh}(x)$ as follows.
\begin{align*}
F^{t}_{nh}(x) &=t F^{-}_{nh}(x) +(1-t) F^+_{nh}(x)
\nonumber\\
&= \sum_{i=1}^n \frac{1}{q(x)}\frac{1}{nh}\Big(t w\Big(\frac{x-V_i}{h}\Big)
-(1-t) w\Big(\frac{x+1-V_i}{h}\Big)\Big)+1-t\\
&= \frac{1}{n}\sum_{i=1}^nV_{ih}(x)+1-t \nonumber
\end{align*}
where
$$
V_{ih}(x)=\frac{1}{q(x)}\frac{1}{h}\Big(t w\Big(\frac{x-V_i}{h}\Big)
-(1-t) w\Big(\frac{x+1-V_i}{h}\Big)\Big).
$$
We are interested in the even moments of $V_{ih}$. The equivalent of Lemma \ref{lem:um} is Lemma \ref{lem:vm} below.
\begin{lem}\label{lem:vm} For $h \to 0$ and $m$ even we have
\begin{equation}\nonumber
\ex V_{ih}(x)^m=\frac{1}{q(x)^{m-1}}\frac{1}{h^{m-1}}(t^m F(x) +(-1)^m(1-t)^m(1-F(x))\int
w(v)^mdv+ o\Big( \frac{1}{h^{m-1}}\Big).
\end{equation}
\end{lem}
The proof is similar as the proof of Lemma \ref{lem:um} and is given in Section \ref{lemmaproof3}.\\

Note that we have 
\begin{equation}
\var F^{t}_{nh}(x) = \var \frac{1}{n}\sum_{i=1}^n\,V_{ih}(x) = \frac{1}{n^2} \var \sum_{i=1}^n\,V_{ih}(x) = \frac{1}{n} \var V_{1h}(x).
\end{equation}
Hence by Lemma \ref{lem:vm} the expansion of the variance of $F^{t}_{nh}(x)$ is
\begin{align*}
\var F^{t}_{nh}(x) &= \frac{1}{n} \var V_{1h}(x)\\
&= \frac{1}{n} \left(\ex V_{ih}(x)^2 - \big(\ex V_{ih}(x)\big)^2 \right)\\
&= \frac{1}{n}\left(\frac{1}{q(x)}\frac{1}{h}(t^2 F(x) +(1-t)^m(1-F(x))\int
w(v)^mdv+ o\Big( \frac{1}{h}\Big) - O(1)^2\right)\\
&= \frac{1}{q(x)}\frac{1}{nh}(t^2 F(x) +(1-t)^2(1-F(x))\int
w(v)^2dv+ o\Big( \frac{1}{nh}\Big).
\end{align*}
$\hfill \Box$

\subsection{Proof of Lemma \ref{lem:vm}}\label{lemmaproof3}
In line with the proof of Lemma \ref{lem:um} in Section \ref{asknownq} one has to define
$$a = t\left(\frac{1}{q(x)} \frac{1}{h}\, w\Big(\frac{x-V_i}{h}\Big)\right)$$
and
$$b=-(1-t)\left(-\frac{1}{q(x)} \frac{1}{h}\, w\Big(\frac{x+1-V_i}{h}\Big)\right).$$\\
Following the same steps we arrive at
\begin{equation}\label{MomentsZ}
\begin{split}
\ex V_{ih}(x)^m &= t^m\frac{1}{q(x)^m}\frac{1}{h^m}\ex [ w^m\Big(\frac{x-V_i}{h}\Big)]\\
& \quad + (1-t)^m\frac{1}{q(x)^m}\frac{1}{h^m}\ex [ w^m\Big(\frac{x + 1 -V_i}{h}\Big)]\\
&= t^m\frac{1}{q(x)^{m-1}}\frac{1}{h^{m-1}}F(x) \int w^m(u)du + o\left(\frac{1}{h^{m-1}}\right)\\
& \quad + (1-t)^m\frac{1}{q(x)^{m-1}}\frac{1}{h^{m-1}}(1-F(x)) \int  w^m(u)du + o\left(\frac{1}{h^{m-1}}\right).
\end{split}
\end{equation}
The last line is proven in the same way as equation (1) and (2) in Section \ref{asknownq} are proven. Lemma \ref{lem:vm} is proven by the last line of equation (\ref{MomentsZ}) . $\hfill \Box$\\

\subsection{Proof of Theorem \ref{ThmUnkownq}}

Recall
\begin{equation}
f_{nh} (x )=\frac{r_{nh}(x)}{q(x)^2},
\end{equation}
where
\begin{align*}
r_{nh}(x)&=q(x)\Big((1-\hat F_n(x))g_{nh}'(x)-\hat F_n(x)g_{nh}'(x+1)\Big)\\
&\quad-q'(x)\Big((1-\hat F_n(x))g_{nh}(x)-\hat F_n(x)g_{nh}(x+1)\Big).
\end{align*}
If we substitute the estimator $q_{n\tilde h}(x)$ for $q(x)$ then we get
\begin{equation}
f_{nh\tilde h} (x )=\frac{r_{nh\tilde h}(x)}{q_{n\tilde h}(x )^2},
\end{equation}
where
\begin{align*}
r_{nh\tilde h}(x)&=q_{n\tilde h}(x)\Big((1-\hat F_n(x))g_{nh}'(x)-\hat F_n(x)g_{nh}'(x+1)\Big)\\
&\quad-q_{n\tilde h}'(x)\Big((1-\hat F_n(x))g_{nh}(x)-\hat F_n(x)g_{nh}(x+1)\Big).
\end{align*}

The first step is a linearisation similar to the linearisation of the Nadaraya Watson estimator in Ha\"{e}rdle (1990), p. 99.
We have
\begin{equation}\label{lin}
f_{nh\tilde h} (x) - f(x)=   \frac{r_{nh\tilde h}(x)- f(x)q_{n\tilde h}(x)^2}{{q (x)^2}}
 +\Big(f_{nh\tilde h} (x)-f(x)\Big)\Big(\frac{q (x)^2-q_{n\tilde h}(x)^2}{q (x)^2}\Big).
\end{equation}
We rewrite the first term as
\begin{align*}\label{t1}
&\frac{r_{nh}(x)+r_{nh\tilde h}(x)-r_{nh}(x)- f(x)q(x)^2
-f(x)(q_{n\tilde h}(x)^2-q(x)^2)}{q(x)^2}\\
&\quad\quad = f_{nh}(x)-f(x)
 +\frac{1}{q(x)^2}\,(r_{nh\tilde h}(x)-r_{nh}(x))
-\frac{f(x)}{q(x)^2}\,(q_{n\tilde h}(x)^2-q(x)^2).
\end{align*}
By the weak consistency of $\hat F_n(x)$, which follows from (\ref{mseconsistency}), $g_{nh}(x)$ and $g_{nh}'(x)$, the difference of $r_{nh\tilde h}(x)$ and $r_{nh}(x)$ can be rewritten as
\begin{align*}
r_{nh\tilde h}(x)&-r_{nh}(x) \\
&=(q_{n\tilde h}(x)-q(x)\Big((1-\hat F_n(x))g_{nh}'(x)-\hat F_n(x)g_{nh}'(x+1)\Big)\\
&-(q_{n\tilde h}'(x)-q'(x))\Big((1-\hat F_n(x))g_{nh}(x)-\hat F_n(x)g_{nh}(x+1)\Big)\\
&=(q_{n\tilde h}(x)-q(x))(f(x)q(x)+o_P(1))-(q_{n\tilde h}'(x)-q'(x))o_P(1).
\end{align*}
We also have
\begin{align*}
&-\frac{f(x)}{q(x)^2}\,(q_{n\tilde h}(x)^2-q(x)^2)=-\frac{f(x)}{q(x)^2}\,\Big((q(x)+q_{n\tilde h}(x)-q(x))^2-q(x)^2\Big)\\ &=-2\frac{f(x)}{q(x)}\,(q_{n\tilde h}(x)-q(x))-\frac{f(x)}{q(x)^2}\,(q_{n\tilde h}(x)-q(x))^2.
&
\end{align*}
Now by
\begin{align*}
&q_{n\tilde h}(x)-q(x)=\frac{1}{2}\,q(x)''{\tilde h}^2\int_{-\infty}^\infty v^2w(v)dv+o({\tilde h}^2)+O_P\Big(\frac{1}{(n{\tilde h})^{1/2}}\Big),\\
&q_{n\tilde h}(x)'-q(x)'=O({\tilde h}^2)+O_P\Big(\frac{1}{(n{\tilde h}^3)^{1/2}}\Big)
\end{align*}
we get
\begin{equation}\label{term1}
\frac{1}{q(x)^2}\,(r_{nh\tilde h}(x)-r_{nh}(x)) =\frac{1}{2}{\tilde h}^2\frac{f(x)q''(x)}{q(x)}\int_{-\infty}^\infty v^2w(v)dv+o({\tilde h}^2)+o_P\Big(\frac{1}{(n{\tilde h}^3)^{1/2}}\Big)
\end{equation}
and
\begin{equation}\label{term2}
-2\frac{f(x)}{q(x)}\,(q_{n\tilde h}(x)^2-q(x)^2)=
-{\tilde h}^2\frac{f(x)q''(x)}{q(x)}\int_{-\infty}^\infty v^2w(v)dv+o({\tilde h}^2)+O_P\Big(\frac{1}{(n{\tilde h})^{1/2}}\Big).
\end{equation}
Taking (\ref{term1}) and (\ref{term2}) together we get
$$
-\frac{1}{2}\,{\tilde h}^2\frac{f(x)q''(x)}{q(x)}\int_{-\infty}^\infty v^2 w(v)dv+o({\tilde h}^2)+o_P\Big(\frac{1}{(n{\tilde h^3})^{1/2}}\Big).
$$
This representation shows that the claims of the theorem hold for the first term in (\ref{lin}).

Finally let us consider the second term in (\ref{lin}).
By the weak consistency of $\hat F_n(x)$, $g_{nh}(x)$ and $g_{nh}'(x)$ the estimator $f_{nh\tilde h} (x)$ is a weakly consistent estimator of $f(x)$. Hence
\begin{equation}
\Big(f_{nh\tilde h} (x)-f(x)\Big)\Big(\frac{q (x)^2-q_{n\tilde h}(x)^2}{q (x)^2}\Big)=o_P(1)O_P\Big(\frac{1}{\sqrt{n\tilde h}}\Big)=o_P(n^{-3/7}),
\end{equation}
which renders this term negligible.
\hfill$\Box$

\bigskip

\noindent{\large \bf Acknowledgement}

\medskip

The major part of the paper is based on the master thesis of Graafland at the Korteweg-de Vries Institute for Mathematics of the University of Amsterdam.

\end{document}